\documentclass[12pt]{article}
\usepackage{amsmath, amsthm, amssymb, stmaryrd, fullpage}

\newtheorem{theorem}{Theorem}[subsection]
\numberwithin{equation}{theorem}
\newtheorem{lemma}[theorem]{Lemma}
\newtheorem{conj}[theorem]{Conjecture}
\newtheorem{cor}[theorem]{Corollary}
\newtheorem{prop}[theorem]{Proposition}
\newtheorem{question}[theorem]{Question}
\theoremstyle{definition}
\newtheorem{defn}[theorem]{Definition}
\newtheorem{example}[theorem]{Example}
\newtheorem{convention}[theorem]{Convention}
\newtheorem{remark}[theorem]{Remark}
\newtheorem{hypothesis}[theorem]{Hypothesis}
\newtheorem{notation}[theorem]{Notation}

\def\FF{\mathbb{F}}
\def\NN{\mathbb{N}}
\def\QQ{\mathbb{Q}}
\def\RR{\mathbb{R}}
\def\ZZ{\mathbb{Z}}
\newcommand{\calE}{\mathcal{E}}
\newcommand{\calF}{\mathcal{F}}
\newcommand{\calO}{\mathcal{O}}
\newcommand{\calR}{\mathcal{R}}
\newcommand{\gotho}{\mathfrak{o}}

\newcommand{\del}{\partial}

\def\be{\mathbf{e}}
\def\bv{\mathbf{v}}
\def\dual{\vee}

\DeclareMathOperator{\Frac}{Frac}
\DeclareMathOperator{\Gal}{Gal}
\DeclareMathOperator{\GL}{GL}
\DeclareMathOperator{\inte}{int}

\DeclareMathOperator{\rank}{rank}
\DeclareMathOperator{\rel}{rel}
\DeclareMathOperator{\sep}{sep}

\DeclareMathOperator{\spect}{sp}
\DeclareMathOperator{\Swan}{Swan}
\DeclareMathOperator{\unr}{unr}

\begin{document}

\title{Swan conductors for $p$-adic differential modules, I: A local
construction}
\author{Kiran S. Kedlaya \\ Department of Mathematics \\ Massachusetts
Institute of Technology \\ 77 Massachusetts Avenue \\
Cambridge, MA 02139 \\
\texttt{kedlaya@mit.edu}}
\date{September 18, 2007}

\maketitle

\begin{abstract}
We define a numerical invariant, the differential Swan conductor, for
certain differential modules on a rigid analytic annulus over a 
$p$-adic field. This gives a definition of a conductor for 
$p$-adic Galois representations with finite local monodromy
over an equal characteristic discretely valued field,
which agrees with the usual Swan conductor when the residue
field is perfect. We also establish analogues of some
key properties of the usual Swan conductor, such as integrality
(the Hasse-Arf theorem), and the fact that the graded pieces of the
associated ramification filtration on Galois groups are abelian 
and killed by $p$.
\end{abstract}

\section*{Introduction}

In this paper, we define a numerical invariant, which we call the differential
Swan conductor, for certain differential modules on a rigid analytic
annulus over a complete nonarchimedean field of mixed characteristics.
We then use this definition to define a differential Swan conductor
for $p$-adic Galois representations with finite local monodromy
over an equal characteristic discretely valued field,
whose residue field need not be perfect. The latter will
coincide with the usual Swan conductor in the case of a perfect residue field.

The construction of the differential Swan conductor proceeds by
measuring the failure of convergence of the Taylor isomorphism, or equivalently,
the failure of local horizontal sections for the connection to converge on as large
a disc as possible. This phenomenon distinguishes the study of differential
equations over $p$-adic fields from its classical analogue, and the relationship
with Swan conductors explains the discrepancy in terms of wild ramification in
characteristic $p$. (The analogy between irregularity of connections and wild
ramification has been known for a while, but recent developments have pushed
it further, e.g., construction of a de Rham analogue of local $\epsilon$-factors
\cite{bbe-epsilon}.)

In the case of Galois representations over an equal characteristic discretely
valued field with perfect
residue field, the differential interpretation of the Swan
conductor is known from the work of several authors, including Andr\'e, 
Christol-Mebkhout, Crew, Matsuda, Tsuzuki, et al; see
\cite[\S 5]{me-mono-over} for an overview. The question of extending this
interpretation
to the case of imperfect residue field was first raised by Matsuda
\cite{matsuda-dwork}, who proposed giving a differential interpretation of the
logarithmic conductor of Abbes-Saito \cite{abbes-saito1, abbes-saito2}.
Our point of view is a bit different: we first construct a numerical invariant
from differential considerations, and check that it has good properties.
These include the Hasse-Arf property, i.e., integrality of conductors
(Theorem~\ref{T:Hasse-Arf}), and the fact that the associated ramification
filtration on Galois groups has graded pieces which are elementary abelian
(Theorem~\ref{T:upper numbering}).
Only then do we pose questions about 
reconciling the definition with other constructions; we do not answer any of
these.

In a subsequent paper, we will
apply this construction to overconvergent $F$-isocrystals
on varieties over perfect fields of positive characteristic;
in particular, the construction applies to discrete representations of
the \'etale fundamental groups of open varieties.
We will
pay particular attention to how the differential Swan conductor of a fixed
isocrystal changes as we vary the choice of a boundary divisor along which
to compute the conductor.

\subsection*{Acknowledgments}
Thanks to Liang Xiao for comments on an early draft.
This material was first presented at the Hodge Theory conference at
Venice International University in June 2006; that presentation 
was sponsored by the Clay Mathematics Institute.
The author was additionally supported by NSF grant DMS-0400727,
NSF CAREER grant DMS-0545904, and a Sloan Research Fellowship.

\section{Differential fields}
\label{sec:diff fields}

We start with a summary of some relevant facts about differential fields
and modules. We defer to \cite[\S 3]{me-part3} (and other
explicitly cited references) for more details.

\subsection{Differential modules and twisted polynomials}

\begin{hypothesis}
Throughout this subsection, 
let $F$ be a differential field of order $1$ and characteristic zero, i.e.,
a field of characteristic zero equipped with a derivation $\del$.
\end{hypothesis}

\begin{defn}
Let $F\{T\}$ denote the (noncommutative)
\emph{ring of twisted polynomials} over $F$
\cite{ore};
its elements are finite formal sums $\sum_{i \geq 0}
 a_i T^i$ with $a_i \in F$, multiplied according to
the rule $Ta = aT + \del(a)$ for $a \in F$.
\end{defn}

\begin{remark} \label{R:opposite}
The opposite ring of $F\{T\}$ is the ring of twisted polynomials for
the differential field given by equipping $F$ with the derivation $-\del$
instead of $\del$.
\end{remark}

\begin{defn}
A \emph{differential module} over $F$ is a finite dimensional $F$-vector
space $V$ equipped with an action of $\del$ (subject to the Leibniz rule);
any such module inherits a left action of $F\{T\}$ where $T$ acts via $\del$.
For $V$ a differential module over $F$,
a \emph{cyclic vector} in $V$ is a vector $\bv \in V$
such that $\bv, \del (\bv), \dots, \del^{\dim(V)-1}
(\bv)$ form a basis of $V$. A cyclic vector defines an isomorphism
$V \cong F\{T\}/F\{T\}P$ of differential modules 
for some twisted polynomial $P \in F\{T\}$, where the $\del$-action
on $F\{T\}/F\{T\}P$ is left multiplication by $T$.
\end{defn}

\begin{lemma} \label{L:cyclic vector}
Every differential module over $F$ contains a 
cyclic vector.
\end{lemma}
\begin{proof}
See, e.g., \cite[Theorem~III.4.2]{g-functions}.
\end{proof}

\begin{hypothesis}
For the remainder of this subsection, assume that the differential
field $F$ is equipped with a nonarchimedean norm $|\cdot|$, and let
$V$ denote a nonzero differential module over $F$.
Write $v(x) = -\log |x|$ for the valuation corresponding to $|\cdot|$.
\end{hypothesis}

\begin{defn}
Let $|\del|_F$ denote the operator norm of $\del$ on $F$.
Let $|\del|_{V,\spect}$ denote the
\emph{spectral norm} of $\del$ on $V$,
i.e., the limit $\lim_{s \to \infty} |\del^s|_{V}^{1/s}$
for any fixed $F$-compatible norm $|\cdot|_V$ on $V$. Any two such norms on 
$V$ are equivalent \cite[Proposition~4.13]{schneider}, 
so the spectral norm does not depend on the choice. More explicitly,
if one chooses a basis of $V$ and lets $D_s$ denote the matrix
via which $\del^s$ acts on this basis, then
\begin{equation} \label{eq:diff module basis}
\max\{|\del|_{F,\spect}, |\del|_{V,\spect} \} = \max\{|\del|_{F,\spect}, 
\lim_{s \to \infty} |D_s|^{1/s} \},
\end{equation}
where the norm applied to $D_s$ is the supremum over entries
\cite[Proposition~1.3]{christol-dwork}.
\end{defn}

\begin{defn}
For $P(T) = \sum_i a_i T^i \in F\{T\}$ a nonzero twisted polynomial, 
define the \emph{Newton polygon} of $P$
as the lower convex hull of the set $\{(-i, v(a_i))\} \subset \RR^2$.
This Newton polygon obeys the usual additivity rules \emph{only} for
slopes less than $-\log |\del|_F$ \cite[Lemma~3.1.5 and
Corollary~3.1.6]{me-part3}, \cite[\S 1]{robba}.
\end{defn}

\begin{prop}[Christol-Dwork] \label{P:read off norm}
Suppose $V \cong F\{T\}/F\{T\}P$ and $P$ has least slope $r$. Then
\[
\max\{|\del|_F, |\del|_{V,\spect}\} = \max\{|\del|_F, e^{-r}\}.
\]
\end{prop}
\begin{proof}
See \cite[Th\'eor\`eme~1.5]{christol-dwork}
or \cite[Proposition~3.3.7]{me-part3}.
\end{proof}

\begin{prop}[Robba] \label{P:slope factor}
Suppose that $F$ is complete for its norm. Then
any monic twisted polynomial $P \in F\{T\}$ admits a unique factorization
\[
P = P_+ P_{m} \cdots P_{1}
\]
such that for some $r_1 < \cdots < r_m < -\log |\del|_F$,
each $P_{i}$ is monic with all slopes equal to $r_i$,
and $P_+$ is monic with all slopes at least $-\log |\del|_F$.
\end{prop}
\begin{proof}
This follows by repeated application of Hensel's lemma for
twisted polynomials \cite{robba}; see also \cite[Corollary~3.2.4]{me-part3}.
\end{proof}

\subsection{Differential fields of higher order}

\begin{hypothesis} \label{H:diff field}
Throughout this subsection, let $F$ denote a \emph{differential field
of order $n$}, i.e., a field $F$ equipped with $n$ commuting
derivations $\del_1, \dots, \del_n$. Assume also that
$F$ has characteristic zero and is complete for a nonarchimedean
norm $|\cdot|$ with corresponding valuation $v$.
Let $V$ denote a nonzero \emph{differential module} over $F$, i.e.,
a nonzero finite dimensional
$F$-vector space  
equipped with commuting actions of $\del_1, \dots, \del_n$.
We apply the results of the previous subsection by
singling out one of $\del_1, \dots, \del_n$.
\end{hypothesis}

\begin{defn}
Define the \emph{scale} of $V$ as
\[
\max\left\{ \max\left\{1, \frac{|\del_i|_{V, \spect}}{|\del_i|_{F, \spect}}
\right\}: 
i \in \{1, \dots, n\} \right\};
\]
note that this quantity is at least 1 by definition, 
with equality at least when $V = F$.
For $i=1, \dots, n$, we say $\del_i$ is \emph{dominant for $V$} if
$\max\{1, |\del_i|_{V, \spect}/|\del_i|_{F, \spect}\}$ equals the scale of $V$.
\end{defn}

\begin{defn}
Let $V_1, \dots, V_m$ be the Jordan-H\"older factors of $V$ (listed
with multiplicity). Define the 
\emph{scale multiset} of $V$ as the multiset of cardinality $\dim_F V$,
consisting of the scale of $V_j$ included
with multiplicity $\dim_F V_j$, for $j=1, \dots, m$.
Note that the largest element of the scale multiset equals the scale of $V$.
\end{defn}

\begin{remark}
If $n=1$ and $V \cong F\{T\}/F\{T\}P$ for $P$ a twisted polynomial, then
Proposition~\ref{P:slope factor} implies that the multiplicity
of any $r< - \log|\del|_F$ as a
slope of the Newton polygon of $P$ coincides with the multiplicity
of $e^{-r}/|\del|_{F,\spect}$ in the scale multiset of $V$.
\end{remark}

\begin{prop} \label{P:scale decomp}
Suppose that $|\del_i|_F/|\del_i|_{F,\spect} = s_0$ for $i=1, \dots, n$.
Then there is a unique decomposition
\[
V = V_- \oplus \bigoplus_{s > s_0} V_s
\]
of differential modules,
such that each Jordan-H\"older factor of $V_s$ has scale $s$,
and each Jordan-H\"older factor of $V_-$ has scale at most $s_0$.
\end{prop}
\begin{proof}
This may be deduced from Proposition~\ref{P:slope factor}, as in 
\cite[Proposition~3.4.3]{me-part3}.
\end{proof}

\begin{defn}
We refer to the decomposition given in Proposition~\ref{P:scale decomp}
as the \emph{scale decomposition} of $V$.
\end{defn}

\section{Conductors for $\nabla$-modules}

In this section, we construct the differential Swan conductor
for certain differential modules over $p$-adic fields.
We will perform all of the calculations under
the bifurcated Hypothesis~\ref{H:Gauss norm1};
one of the two options therein
allows for nonarchimedean fields which are
not discretely valued, but restricts their residue fields,
while the other
is less restrictive on residue fields, but requires the nonarchimedean
norms to be discretely valued.

\setcounter{theorem}{0}
\begin{notation}
For $S$ a set or multiset, write $S^p = \{s^p: s \in S\}$.
If $A,B$ are two multisets of the same cardinality $d$, then write
$A \geq B$ to mean that for $i=1, \dots, d$, the $i$-th largest element
of $A$ is greater than or equal to the $i$-th largest element of $B$
(counting multiplicity).
\end{notation}

\subsection{Setup}

\begin{defn}
Given a field $K$ equipped with a (possibly trivial) nonarchimedean norm, 
for $\rho_1, \dots, \rho_n \in (0, +\infty)$, the
\emph{$(\rho_1, \dots, \rho_n)$-Gauss norm} on $K[u_1, \dots, u_n]$
is the norm $|\cdot|_\rho$ given by
\[
\left| \sum_I c_I u_1^{i_1} \cdots u_n^{i_n} \right|
= \max_I \{|c_I| \rho_1^{i_1} \cdots \rho_n^{i_n} \};
\]
this norm extends uniquely to $K(u_1, \dots, u_n)$.
\end{defn}

\begin{defn}
For $\ell/k$ an extension of fields of characteristic $p>0$, a
\emph{$p$-basis} of $\ell$ over $k$ is a set $B \subset \ell$ with the property that the
products $\prod_{b \in B} b^{e_b}$, where $e_b \in \{0, \dots, p-1\}$
for all $b \in B$ and $e_b = 0$ for all but finitely many $b$, are all distinct
and form a basis for $\ell$ as a vector space over the compositum
$k \ell^p$. 
By a \emph{$p$-basis} of $\ell$, we mean a $p$-basis of $\ell$ over 
$\ell^p$.
\end{defn}

\begin{hypothesis} \label{H:Gauss norm1}
For the rest of this section, assume
\emph{one} of the following two sets of hypotheses.
\begin{enumerate}
\item[(a)]
Let $K$ be a field of characteristic zero, complete for a
(not necessarily discrete)
nonarchimedean norm $|\cdot|$, with residue field $k$
of characteristic $p>0$.
Equip $K(u_1,\dots,u_n)$ with the $(1,\dots,1)$-Gauss norm.
Let $\ell$ be a finite separable extension of $k(u_1,\dots,u_n)$,
and let $L$ be the unramified extension with residue field $\ell$
of the completion of $K(u_1,\dots,u_n)$.
\item[(b)]
Let $K$ be a field of characteristic zero, complete for a
nonarchimedean norm $|\cdot|$, with discrete value
group and residue field $k$ of characteristic $p>0$.
Let $L$ be an extension of $K$, complete for an extension of $|\cdot|$
with the same value group, whose residue field $\ell$ admits
a finite $p$-basis $\overline{B} = \{\overline{u_1}, \dots, \overline{u_n}\}$
over $k$. For $i=1,\dots,n$,
let $u_i$ be a lift of $\overline{u_i}$ to the valuation ring $\gotho_L$
of $L$.
\end{enumerate}
\end{hypothesis}

\begin{defn}
Under either option in Hypothesis~\ref{H:Gauss norm1},
the module of continuous differentials $\Omega^1_{L/K}$ is generated
by $du_1,\dots,du_n$; let $\del_1, \dots, \del_n$ denote the
dual basis of derivations (that is, $\del_i = \frac{\del}{\del u_i}$).
\end{defn}

\begin{remark}
Note that $|\del_i|_L/|\del_i|_{L,\spect} = |p|^{-1/(p-1)}$ 
for $i=1, \dots, n$,
so Proposition~\ref{P:scale decomp} applies.
\end{remark}

\subsection{Taylor isomorphisms}

The scale of a differential module over $L$ can be interpreted as a normalized
radius of convergence for the Taylor series, as follows.

\begin{convention}
Let $\NN_0$ denote the monoid of nonnegative integers.
For $I \in \NN_0^n$ 
and $*$ any symbol, we will write $*^I$ as shorthand for
$*_1^{i_1} \cdots *_n^{i_n}$. We also write $I!$ as shorthand for
$i_1! \cdots i_n!$.
\end{convention}

\begin{defn} \label{D:Taylor inter}
Let $V$ be a differential module over $L$.
Define the \emph{formal Taylor isomorphism} on $V$ to be the
map $T: V \mapsto V \otimes_L L \llbracket x_1, \dots, x_n \rrbracket$ given by
\[
T(\bv) = \sum_{I \in \NN_0^n} \frac{x^I}{I!} \del^I(\bv).
\]
We can then interpret the scale of $V$ as the minimum $\lambda$ such that
$T$ takes values in $V \otimes_L R$, for $R$ the subring of
$L \llbracket x_1, \dots, x_n \rrbracket$ consisting of series
convergent on the open polydisc
\[
|x_i| < \lambda^{-1} \qquad (i=1, \dots, n).
\]
In particular, if $L'$ is a complete extension of $L$, and
$x_1, \dots, x_n \in L'$ satisfy $|x_i| < \lambda^{-1}$ for $\lambda$
the scale of $V$, we obtain a \emph{concrete Taylor isomorphism}
$T(\bv; x_1, \dots, x_n): V \to V \otimes_L L'$ by substitution.
\end{defn}

\begin{remark}
If $x_1, \dots, x_n \in L$ satisfy $|x_i| < 1$, then the concrete
Taylor isomorphism $T(\cdot; x_1, \dots,x_n)$ is defined on $L$, and 
is a $K$-algebra homomorphism carrying $u_i$ to $u_i + x_i$. 
If $V$ is a differential module
of scale $\lambda$, and $|x_i| < \lambda^{-1}$ for $i=1,\dots,n$,
then the concrete Taylor isomorphism $T(\cdot; x_1, \dots, x_n)$ on $V$
is semilinear over the concrete Taylor isomorphism on $L$.
\end{remark}

\begin{remark} \label{R:same scale}
Note that for any $I \in \NN_0^n$, $|\del^I/I!|_F \leq 1$. Hence
if $x_1, \dots, x_n \in L$ satisfy $|x_i| < 1$, then for any
$f \in L$, 
\[
|T(f; x_1, \dots, x_n) - f| \leq \max_i \{|x_i|\} \cdot |f|.
\]
In particular, suppose $u'_1, \dots, u'_n \in L$
satisfy $|u'_i| = 1$, and the images of $u'_1, \dots, u'_n$ in $\ell$
form a $p$-basis of $\ell$ over $k$.
Then $T(\cdot; x_1, \dots, x_n)$ can also be interpreted as the
concrete Taylor isomorphism defined with respect to the dual basis
of $du'_1, \dots, du'_n$ and evaluated at $y_1, \dots, y_n$, for
$y_i = T(u'_i; x_1, \dots, x_n) - u'_i$. 
This implies that the scale of a differential module computed with respect to
(the dual basis to) $du_1, \dots, du_n$ is no greater than with respect to
$du'_1, \dots,du'_n$; by the same calculation in reverse, it follows
that the two scales are equal. (Francesco Baldassarri has suggested a
coordinate-free definition of the scale that explains this remark; we will
follow up on this suggestion elsewhere.)
\end{remark}

\subsection{Frobenius descent}

As discovered originally by Christol-Dwork \cite{christol-dwork},
in the situations of Hypothesis~\ref{H:Gauss norm1},
one can overcome the limitation on scales imposed by 
Proposition~\ref{P:read off norm} by using descent along 
the substitution $u_i \mapsto u_i^p$.

\begin{defn} \label{D:frob ant}
Let $V$ be a differential module over $L$ with scale less than 
$|p|^{-1/(p-1)}$.
If $K$ contains a primitive $p$-th root of unity $\zeta$, we may
define an action of the group $(\ZZ/p\ZZ)^n$ on $V$ using
concrete Taylor isomorphisms:
\[
\bv^J = T(\bv; (\zeta^{j_1} - 1) u_1, \dots, (\zeta^{j_n} - 1) u_n)
\qquad (J \in (\ZZ/p\ZZ)^n).
\]
Let $V_1$ be the fixed space under this group action; in particular,
taking $V = L$, we obtain a subfield $L_1$ of $L$, which we may view
as a differential field of order $n$ for the derivations $\del_{i,1} 
= \frac{\del}{\del (u_i^p)}$. In general,
$V_1$ may be viewed as a differential module over $L_1$, the natural map
$V_1 \otimes_{L_1} L \to V$ is an isomorphism of $L$-vector spaces
(by Hilbert 90), 
and the actions of
$\del_i$ and $\del_{i,1}$ on $V$ are related by the formula
\begin{equation} \label{eq:two deltas}
\del_{i,1} = \frac{1}{p u_i^{p-1}} \del_i.
\end{equation}
We call $V_1$ the \emph{Frobenius antecedent} of $V$.
If $K$ does not contain a primitive $p$-th root of unity,
we may still define the Frobenius antecedent using Galois descent.
\end{defn}

\begin{prop} \label{P:frob scale}
Let $V$ be a differential module over $L$ with scale $s < |p|^{-1/(p-1)}$
and scale multiset $S$. 
Then the scale multiset of the Frobenius antecedent of $V$ is $S^p$.
\end{prop}
\begin{proof}
Since any direct sum decomposition commutes with the formation of the Frobenius
antecedent $V_1$, it suffices to check that the scale of $V_1$ is $s^p$.
Let $T(\bv)$ be the formal Taylor isomorphism for $V$, and let
$T'(\bv)$ be the formal Taylor isomorphism for $V_1$ but with variables
$x'_1, \dots, x'_n$.

By \cite[Lemma~5.12]{me-mono-over},
for $t,t_1$ in any nonarchimedean field,
\begin{equation} \label{eq:p-power}
|t - t_1| < \lambda^{-1} |t| \implies
|t^p - t_1^p| < \lambda^{-p} |t|^p
\qquad (1 < \lambda < |p|^{-1/(p-1)}).
\end{equation}
(We repeat from \cite[Lemma~4.4.2]{me-part3} the description of a 
misprint in the last line of the statement of
\cite[Lemma~5.12]{me-mono-over}: one must read $r^{1/p} \rho^{1/p}, r \rho$
for $r \rho^{1/p}, r^p \rho$, respectively.), 
Hence the convergence of $T'(\bv; x'_1, \dots, x'_n)$
for $|x'_i| < \lambda^{-p}$ implies convergence of
$T(\bv; x_1, \dots, x_n)$ for $|x_i| < \lambda^{-1}$,
so the scale of
$V_1$ is at least $s^p$. On the other hand, we can obtain 
$T'$ by averaging $T$ over the action of $(\ZZ/p\ZZ)^n$,
so the scale of $V_1$ is at most $s^p$.
(Compare \cite[Theorem~6.15]{me-mono-over}.)
\end{proof}

\begin{remark}
It should also be possible to prove Proposition~\ref{P:frob scale}
by raising both sides of \eqref{eq:two deltas} to a large power and
comparing the results,
but this would appear to be somewhat messy.
\end{remark}

\begin{defn}
If $V$ is a differential module over $ F$ of scale less than
$|p|^{-1/(p^{m-1}(p-1))}$, by Proposition~\ref{P:frob scale},
we can iterate the construction of a Frobenius antecedent $m$ times; we call
the result the \emph{$m$-fold Frobenius antecedent} of $V$.
\end{defn}

\begin{remark}
Note that it is also possible to construct antecedents
one variable at a time; the point is that since the operators $\del_i$, 
$\del_j$ commute for $i \neq j$, $\del_i$ continues to act on the
antecedent with respect to $\del_j$. This will be used in the
proof of Proposition~\ref{P:break piecewise}.
\end{remark}

\subsection{$\nabla$-modules} 

\begin{notation}
Let $\Gamma^*$ denote the divisible closure of $|K^*|$. We say a
subinterval of $(0, +\infty)$ is \emph{aligned} if each endpoint at which
it is closed belongs to $\Gamma^*$.
\end{notation}

\begin{remark}
One can drop the word ``aligned'', and all references to $\Gamma^*$,
everywhere hereafter if one works with Berkovich analytic spaces 
\cite{berkovich} instead of rigid
analytic spaces. We omit further details.
\end{remark}

\begin{notation}
For $I$ an aligned interval and $t$ a dummy variable, 
let $A_{L}(I)$ be the rigid analytic (over $L$) subspace of the affine $t$-line 
over $L$ consisting
of points with $|t| \in I$; this space is affinoid if $I$ is closed.
(We omit the parentheses if $I$ is described explicitly, e.g.,
if $I = [\alpha, \beta)$, we write $A_L[\alpha,\beta)$ for $A_L(I)$.)
For $\rho \in I$, we write $|\cdot|_\rho$ for the $\rho$-Gauss norm
\[
\left| \sum_{i \in \ZZ} c_i t^i \right|_\rho = \sup_i \{|c_i| \rho^i\};
\]
for $\rho \in \Gamma^*$, we may interpret $|\cdot|_\rho$ as the supremum
norm on the affinoid space $A_L[\rho,\rho]$.
\end{notation}

\begin{lemma} \label{L:three circles}
Let $I$ be an aligned interval.
For $\rho, \sigma \in I$ 
and $c \in [0,1]$, put $\tau = \rho^c \sigma^{1-c}$. Then
for any $f \in \Gamma(A_L(I), \calO)$,
\[
|f|_\tau \leq |f|_\rho^c |f|_\sigma^{1-c}.
\]
\end{lemma}
\begin{proof}
See \cite[Lemma~3.1.6]{me-part1},
\cite[Corollaire~4.2.8]{amice}, or
\cite[Corollaire~5.4.9]{christol-robba}.
\end{proof}

\begin{defn}
For $I$ an aligned interval, a \emph{$\nabla$-module} on $A_L(I)$
(relative to $K$) is a coherent locally free sheaf $\calE$ on $A_L(I)$
equipped with an integrable $K$-linear connection $\nabla: \calE \to
\calE \otimes \Omega^1_{A_L(I)/K}$. (Here $\Omega^1_{A_L(I)/K}$ denotes the sheaf of
continuous differentials; it is freely
generated over $\calO_{A_L(I)}$
by $du_1, \dots, du_n, dt$.) The connection equips $\calE$
with actions of the derivations $\del_i = \frac{\del}{\del u_i}$ for
$i=1,\dots,n$ and $\del_{n+1} = \frac{\del}{\del t}$; integrability
of the connection is equivalent to commutativity between these actions.
\end{defn}

\begin{defn} \label{D:gen radius}
For $I$ an aligned interval
and $\rho \in I$, let $F_\rho$ be the completion
of $L(t)$ for the $\rho$-Gauss norm, viewed as a differential field
of order $n+1$. For $\calE$ a nonzero $\nabla$-module on $A_L(I)$,
let $J$ be a closed aligned neighborhood of $\rho$ in $I$,
and put
\[
\calE_\rho = \Gamma(A_L(J), \calE)
\otimes_{\Gamma(A_L(J), \calO)} F_\rho,
\]
viewed as a differential module over $F_\rho$; this construction
does not depend on $J$.
Define the \emph{radius multiset} of $\calE_\rho$, denoted
$S(\calE,\rho)$, as the multiset of reciprocals of the scale multiset
of $\calE_\rho$.
Define the \emph{(toric) generic radius of convergence} of $\calE_\rho$,
denoted $T(\calE,\rho)$, as the smallest element of $S(\calE,\rho)$, i.e.,
the reciprocal of the scale of $\calE_\rho$.
\end{defn}

\begin{remark}
As in \cite{me-part3}, the toric generic radius of convergence is normalized
differently from the generic radius of convergence of Christol-Dwork
\cite{christol-dwork}, which would be multiplied by an extra factor of $\rho$.
Our chief justification for this normalization is ``because it works'', in the sense of
giving the expected answer for Example~\ref{ex:artin-schreier}. We
look forward to ongoing work of Baldassarri (compare Remark~\ref{R:same scale})
for a more intrinsic justification.
\end{remark}

\begin{remark}
To our knowledge, 
the consideration of $\nabla$-modules over a rigid analytic
 annulus, but taking into account
derivations of the base field over a subfield, is novel to
this paper. It may prove an interesting exercise to transcribe
the arguments of \cite{me-mono-over}, such as local duality,
 as much as possible to this setting.
\end{remark}

\subsection{The highest ramification break}

\begin{defn}
Let $\calE$ be a $\nabla$-module on $A_L(\epsilon,1)$ for some
$\epsilon \in (0,1)$. We say
$\calE$ is \emph{solvable at $1$} if
\[
\lim_{\rho \to 1^-} T(\calE,\rho) = 1.
\]
\end{defn}

\begin{hypothesis} \label{H:solvable}
For the rest of this subsection,
let $\calE$ be a $\nabla$-module on $A_L(\epsilon,1)$ for some
$\epsilon \in (0,1)$, which is solvable at $1$.
\end{hypothesis}

\begin{lemma} \label{L:concave}
For each $i \in \{1,\dots,n+1\}$,
for $r \in (0, -\log \epsilon)$, put $\rho = e^{-r}$ and
let $f_i(r)$ be the negative logarithm of the scale of $\del_i$
on $\calE_\rho$.
Then $f_i$ is a concave function of $r$; in particular,
\[
\log T(\calE, e^{-r}) = \min_i \{f_i(r)\}
\]
is a concave function of $r$. (This does not require solvability at $1$.)
\end{lemma}
\begin{proof}
It suffices to check concavity on $- \log(J)$ for $J$ an
arbitrary closed aligned subinterval of $(\epsilon,1)$.
Since $J$ is closed aligned, $A_L(J)$ is affinoid; by Kiehl's theorem
(see for instance \cite[Theorem~4.5.2]{fvdp}),
$\Gamma(A_L(J), \calE)$ is a finitely generated module over the ring
$\Gamma(A_L(J), \calO)$. Since that ring is a principal ideal domain
\cite[Proposition~4, Corollaire]{lazard}, 
$\Gamma(A_L(J), \calE)$ is freely generated by some subset
$\be_1, \dots, \be_m$.
Let $D_{i,l}$ be the matrix over $\Gamma(A_L(J), \calO)$ via which
$\del_i^l$ acts on $\be_1, \dots, \be_m$.
For $\rho, \sigma \in J$ and $c \in [0,1]$, put
$\tau = \rho^c \sigma^{1-c}$.
By Lemma~\ref{L:three circles}, we have
\[
|D_{i,l}|_\tau \leq |D_{i,l}|^c_\rho
|D_{i,l}|^{1-c}_\sigma;
\]
taking $l$-th roots of both sides and taking limits yields
\[
\lim_{l \to \infty} |D_{i,l}|_\tau^{1/l} \leq
\left( \lim_{l \to \infty} |D_{i,l}|_\rho^{1/l} \right)^c
\left( \lim_{l \to \infty} |D_{i,l}|_\sigma^{1/l} \right)^{1-c}.
\]
By \eqref{eq:diff module basis}, this yields the desired result.
(Compare \cite[Proposition~4.2.6]{me-part3}.)
\end{proof}

\begin{prop} \label{P:break piecewise}
The function $f(r) = \log T(\calE, e^{-r})$ on $(0, -\log \epsilon)$ is
piecewise linear, with slopes in $(1/(\rank \calE)!) \ZZ$.
Moreover, $f$ is linear in a neighborhood of $0$.
\end{prop}
\begin{proof}
Since $f$ is concave by Lemma~\ref{L:concave}, 
takes nonpositive values, and tends to 0 as
$r \to 0^+$, it is everywhere nonincreasing. Hence
for sufficiently large integers $h$, we can choose $\rho_h \in
(\epsilon,1)$ with $T(\calE,\rho_h) = |p|^{1/(p^{h-1}(p-1))}$ and
$\rho_h < \rho_{h+1}$. Put $r_h = -\log \rho_h$.

We now check piecewise linearity and the slope restriction
on $(r_{h+1},r_h)$; it suffices
to check on $-\log(J)$ for $J$ an
arbitrary closed aligned subinterval of $(\rho_h,\rho_{h+1})$. 
Assume without loss of generality that $K$ contains a primitive $p$-th
root of unity.
Put $L_0 = L$. For $l=1, \dots, h$, let $L_l$ be the subfield of $L_{l-1}$
fixed under the action of $(\ZZ/p\ZZ)^n$ given in 
Definition~\ref{D:frob ant}, but with $u_1^{p^{l-1}}, \dots, u_n^{p^{l-1}}$ 
playing the
roles of $u_1, \dots, u_n$.
Since $T(\calE,\rho) > |p|^{1/(p^{h-1}(p-1))}$ for $\rho \in J$,
using Definition~\ref{D:frob ant} (in the $u_1,\dots,u_n$-directions)
and \cite[Theorem~6.15]{me-mono-over} (in the $t$-direction),
we can construct an $h$-fold Frobenius antecedent $\calE_h$ for
$\calE$, which is defined on $A_{L_h}(J^{p^h})$. 

Apply Lemma~\ref{L:cyclic vector} to construct
a cyclic vector for $\calE_h$ over $\Frac 
\Gamma(A_{L_h}(J^{p^h}), \calO)$;
by writing down the corresponding twisted polynomial $P(T)$ and applying
Proposition~\ref{P:read off norm}, we see that
for $\sigma \in J^{p^h}$, 
$T(\calE_h,\sigma)$ is piecewise of the form
$|g|_{\sigma}^{1/j}$ for some $g \in \Frac \Gamma(A_{L_h}(J^{p^h}), \calO)$
and $j \in \{1, \dots, \rank(\calE)\}$. 
In particular, for $\sigma = \rho^{p^h}$, 
this expression is piecewise of the form
$(|a| \rho^{i p^h})^{1/j}$ for some $a \in K^*$, $i \in \ZZ$,
and $j \in \{1, \dots, \rank(\calE)\}$. 
This proves that on $(r_{h+1}, r_h)$,
$f$ is piecewise linear with slopes in $(1/(\rank \calE)!)\ZZ$.

To check piecewise linearity in a neighborhood of $r_{h}$,
note that as we approach $r_{h}$ from the right, the successive slopes
of $f$ that we encounter are increasing but bounded above, and lie in a
 discrete subset of $\RR$. Hence they stabilize, so $f$ is linear in a
one-sided neighborhood of $r_{h}$. An analogous argument applies again
when approaching
$r_{h+1}$ from the left, so $f$ is piecewise linear
on $[r_{h+1}, r_h]$; taking the union of these intervals,
we deduce that $f$ is piecewise linear on $(0, r_h]$ for some $h$.
An analogous argument applies yet again
when approaching 0 from the right, yielding the desired result.
(Compare \cite[Th\'eor\`eme~4.2-1]{cm3}.)
\end{proof}
\begin{cor} \label{C:highest break}
There exists
$b \in \QQ_{\geq 0}$ 
such that $T(\calE,\rho) = \rho^b$ for all $\rho \in (\epsilon,1)$.
\end{cor}

\begin{defn}
We will refer to the number $b$ in Corollary~\ref{C:highest break} 
as the \emph{(differential) highest ramification break} of $\calE$,
denoted $b(\calE)$.
\end{defn}

\subsection{Invariance}

\begin{defn}
Define the \emph{Robba ring} over $L$ as
\[
\calR_L = \cup_{\epsilon \in (0,1)} \Gamma(A_L(\epsilon,1), \calO).
\]
The elements of $\calR_L$ can be represented as formal Laurent
series $\sum_{i \in \ZZ} c_i t^i$
with $c_i \in L$;
let $\calR_L^{\inte}$ be the subring of series with $|c_i| \leq 1$
for all $i \in \ZZ$.
The ring $\calR_L^{\inte}$ is local, with
maximal ideal consisting of series
with $|c_i| < 1$ for all $i \in \ZZ$, with residue field
$\ell((t))$.
\end{defn}

We first examine invariance under certain endomorphisms of $L$,
following Definition~\ref{D:Taylor inter}.
\begin{defn} \label{D:maps}
Choose $u'_1, \dots, u'_n, t' \in \calR_L^{\inte}$ such that 
under the projection $\calR_L^{\inte} \to \ell((t))$,
$u'_1 - u_1, \dots, u'_n - u_n$ map into $t \ell \llbracket t \rrbracket$
and $t'-t$ maps into $t^2 \ell \llbracket t \rrbracket$.
Then for some $\epsilon \in (0,1)$, the Taylor series
\[
\sum_{I \in \NN^{n+1}_{0}} \frac{(u'_1-u_1)^{i_1}
\cdots (u'_n - u_n)^{i_n} (t'-t)^{i_{n+1}}}{I!}
\del^I (f)
\]
converges for $f \in \Gamma(A_L(I), \calO)$ for any closed aligned
subinterval $I$ of $(\epsilon,1)$, so we can use it to define a map
$g: \Gamma(A_L(I), \calO) \to \Gamma(A_L(I), \calO)$ such that
$g^*(u_i) = u'_i$, $g^*(t) = t'$.
\end{defn}

\begin{prop} \label{P:maps}
Let $g$ be a map as in Definition~\ref{D:maps}.
For any
$\nabla$-module $\calE$ on $A_L(\epsilon,1)$ which is solvable at $1$, 
we have $T(\calE,\rho) = T(g^* \calE, \rho)$ for all
$\rho \in (\epsilon,1)$ sufficiently close to $1$.
In particular, $g^* \calE$ is also solvable at $1$,
and $\calE$ and $g^* \calE$ have the same highest break.
\end{prop}
\begin{proof}
By the choice of $u'_1, \dots, u'_n, t'$, for
$\rho \in (0,1)$ sufficiently close to 1,
\[
|u'_i-u_i|_\rho < 1 \qquad (i=1, \dots, n), \qquad
|t' - t|_\rho < \rho.
\]
We will prove the claim for such $\rho$.

By continuity of $T(\calE,\rho)$ (implied by Lemma~\ref{L:concave}),
it suffices to check for $\rho \in \Gamma^*$. 
There is no loss of generality in enlarging $K$, so we may in fact
assume that there exists $\lambda \in K$ with $|\lambda| = \rho$.
In this case, we may put ourselves in the situation of
Remark~\ref{R:same scale} by 
considering $\calE_\rho$ to be a differential module over the
the completion of $L(t/\rho)$ for the 1-Gauss norm, comparing
the $p$-bases $u_1,\dots,u_n,t/\lambda$ and
$u'_1,\dots,u'_n,t'/\lambda$. This yields the claim.
\end{proof}

\begin{prop} \label{P:maps2}
Let $g: A_L(I) \to A_L(I)$ be the map fixing $L$ and pulling $t$ back to
$t^{p^N}$ for some positive integer $N$.
Then for any
$\nabla$-module $\calE$ on $A_L(\epsilon,1)$, 
we have $S(\calE,\rho) \leq S(g^* \calE, \rho^{1/p^N})$
for all $\rho \in (\epsilon,1)$; moreover, if $n=0$, then
$S(\calE,\rho) \leq S(g^* \calE, \rho^{1/p^N})^{p^N}$.
\end{prop}
\begin{proof}
If we compare the scale multisets of $\del_i$ on $\calE_\rho$ and
on $(g^* \calE)_{\rho^{1/p^N}}$, then we get identical results
for $i=1, \dots, n$. For $i=n+1$, the scale multiset on $\calE_\rho$ is at least
the $p^N$-th power of the scale multiset on $(g^* \calE)_{\rho^{1/p^N}}$,
as in the proof of Proposition~\ref{P:frob scale}.
This yields the claim.
\end{proof}

\begin{prop} \label{P:maps3}
Let $g: A_L(I) \to A_L(I)$ be the map fixing $L$ and pulling $t$ back to
$t^{N}$ for some positive integer $N$ coprime to $p$.
Then for any
$\nabla$-module $\calE$ on $A_L(\epsilon,1)$, 
we have $S(\calE,\rho) = S(g^* \calE, \rho^{1/N})$
for all $\rho \in (\epsilon,1)$.
\end{prop}
\begin{proof}
If we compare the scale multisets of $\del_i$ on $\calE_\rho$ and
on $(g^* \calE)_{\rho^{1/N}}$, then we get identical results
for $i=1, \dots, n$. For $i=n+1$, we again get identical results 
by virtue of \cite[Lemma~5.11]{me-mono-over}.
\end{proof}

We next examine what happens when we change the $p$-basis.
\begin{prop} \label{P:maps1}
Choose $u'_1, \dots, u'_n \in \calR_L^{\inte}$ such that 
under the projection $\calR_L^{\inte} \to \ell((t))$,
$u'_1, \dots, u'_n$ map to elements of $\ell \llbracket t \rrbracket$
lifting a $p$-basis of $\ell$ over $k$.
Let
$\del'_1, \dots, \del'_n$ be the derivations dual to
the basis $du'_1,\dots,du'_n$ of $\Omega^1_{L/K}$.
Let $\calE$ be a $\nabla$-module on $A_L(\epsilon,1)$ for some
$\epsilon \in (0,1)$, which is solvable at $1$.
Then for $\rho \in (0,1)$ sufficiently close to $1$,
the scale of $\calE_\rho$ for $\del_1, \dots, \del_n, \del_{n+1}$
is the same as for $\del'_1,\dots,\del'_n,\del_{n+1}$;
in particular, the highest break is the same in both cases.
\end{prop}
\begin{proof}
If $u'_1, \dots, u'_n \in L$, then we can invoke Remark~\ref{R:same scale}
to obtain the claim. In general, we may first make a transformation as in
the previous sentence, to match up the reductions modulo $t \ell 
\llbracket t \rrbracket$, then invoke Proposition~\ref{P:maps}.
\end{proof}

\subsection{The break decomposition}
\label{subsec:break decomp}

Throughout this subsection, retain Hypothesis~\ref{H:solvable}.

\begin{defn}
We say that $\calE$ has a \emph{uniform break} if
 for all $\rho \in (0,1)$ sufficiently close to 1, $S(\calE,\rho)$ consists
of a single element with multiplicity $\rank(\calE)$. We write
``$\calE$ has uniform break $b$'' as shorthand for
``$\calE$ has a uniform break and its highest ramification break is $b$''.
\end{defn}

\begin{theorem} \label{T:break decomp}
For some $\eta \in (0,1)$, there exists a decomposition of
$\nabla$-modules (necessarily unique)
$\calE = \oplus_{b \in \QQ_{\geq 0}} \calE_b$ over $A_L(\eta,1)$ such that
each $\calE_b$ has uniform break $b$.
\end{theorem}
We will prove Theorem~\ref{T:break decomp} later in this subsection. To begin with,
we recall that the case $L = K$ is essentially a theorem of Christol-Mebkhout
\cite[Corollaire~2.4-1]{cm4},
from which we will bootstrap to the general case.
\begin{lemma} \label{L:cm}
Theorem~\ref{T:break decomp} holds in case $L=K$.
\end{lemma}
\begin{proof}
This is the conclusion of \cite[Corollaire~2.4-1]{cm4}, at least in case $K$
is spherically complete. However, it extends to the general case as follows.

By a straightforward application of Zorn's lemma, we may
embed $K$ into a spherically complete field $K'$.
Apply \cite[Corollaire~2.4-1]{cm4} to obtain a 
break decomposition over $A_{K'}(\eta,1)$ for some $\eta \in (0,1)$;
let $\bv \in \Gamma(A_{K'}(\eta,1), \calE^\dual \otimes \calE)$ be the projector
onto the highest break component.

Now set notation as in the proof of Proposition~\ref{P:break piecewise}.
The set of $\rho \in (\rho_h, \rho_{h+1})$ for which at least one coefficient
$P(T)$ fails to be a unit in $A_{L_h}[\rho^{p^h},\rho^{p^h}]$ 
is discrete, so we may choose $\rho \in (\rho_h, \rho_{h+1})$ not of that form.
Then Proposition~\ref{P:slope factor} gives a factorization of $P(T)$ over
$A_{L_h}[\rho^{p^h},\rho^{p^h}]$ (and likewise in the opposite ring); we thus obtain
an element $\bv'$ of $\Gamma(A_K[\rho,\rho], \calE^\dual \otimes \calE)$
which agrees with $\bv$ over $A_{K'}[\rho,\rho]$. 

For any closed aligned subinterval
$J$ of $(\eta,1)$
containing $\rho$, inside $\Gamma(A_{K'}[\rho,\rho], \calO)$ we have
\[
\Gamma(A_K[\rho,\rho],\calO) \cap \Gamma(A_{K'}(J), \calO) = 
\Gamma(A_K(J), \calO).
\]
Since $\calE^\dual \otimes \calE$ is free over $A_K(J)$ (as in the proof
of Lemma~\ref{L:concave}), this implies that
\[
\Gamma(A_K[\rho,\rho],\calE^\dual \otimes
\calE) \cap \Gamma(A_{K'}(J), \calE^\dual \otimes \calE) = 
\Gamma(A_K(J), \calE^\dual \otimes \calE),
\]
and so $\bv \in \Gamma(A_K(J), \calE^\dual \otimes \calE)$.
Running this argument over all possible $J$, we obtain
$\bv \in \Gamma(A_K(\eta,1), \calE^\dual \otimes \calE)$, so $\calE$
admits a break decomposition over $A_K(\eta,1)$ as desired.
\end{proof}

We exploit Lemma~\ref{L:cm} via the following construction.
\begin{defn}
Define the \emph{relativization} $\calF$ of $\calE$ as the $\nabla$-module
$\calE$ itself, but viewed relative to $L$ instead of $K$. That is, retain only
the action of $\del_{n+1}$. (The term ``generic fibre'' was used in an
earlier version of this paper, but we decided to reserve that name for a
different concept to appear in a subsequent paper.)
\end{defn}

However, we are forced to make a crucial distinction.
\begin{lemma} \label{L:eventually dominant}
For $i \in \{1, \dots, n+1\}$, there exists $\eta \in (0,1)$
such that one of the following two statements is true.
\begin{itemize}
\item For all $\rho \in (\eta,1)$,
$\del_i$ is dominant for $\calE_\rho$.
\item For all $\rho \in (\eta,1)$,
$\del_{i}$ is \emph{not} dominant for $\calE_\rho$.
\end{itemize}
\end{lemma}
\begin{proof}
Let $b$ denote the highest break of $\calE$.
Choose $\eta \in (0,1)$ such that $T(\calE,\rho) = \rho^b$
for all $\rho \in (\eta,1)$.
Put
\[
f_i(\rho) = \frac{|\del_i|_{F_\rho, \spect}}{|\del_i|_{\calE_\rho, \spect}};
\]
then Lemma~\ref{L:concave} shows that
$f_i$ is log-concave. Consequently, if $f_i(\rho) = T(\calE,\rho)$
for two different values of $\rho$, then the same is true for all
intermediate values. This proves the claim: if the second
statement does not hold, then 
there exist $\rho \in (0,1)$ arbitrarily close
to 1 such that $f_i(\rho) = T(\calE,\rho)$, in which case
the first statement holds with $\eta$ equal to any such $\rho$.
\end{proof}

\begin{defn}
For $i \in \{1, \dots, n+1\}$, we say that $\del_i$ is
\emph{eventually dominant} for $\calE$ if the first alternative
in Lemma~\ref{L:eventually dominant} holds, i.e.,
if there exists $\eta \in (0,1)$ such that for all $\rho \in (\eta,1)$,
$\del_i$ is dominant for $\calE_\rho$.
\end{defn}

\begin{remark} \label{R:easy case}
Note that if $\del_{n+1}$ is eventually dominant for $\calE$, then the
highest break term in the decomposition of $\calF$ (which is respected
by $\del_1, \dots, \del_n$ because it is unique)
already has a uniform break. 
Our strategy in case $\del_{n+1}$ is not eventually dominant for
$\calE$ is to perform an operation which one might call \emph{rotation}
to recover that more favorable
situation: namely, we use a concrete Taylor isomorphism to
change the embedding of $K$ into $L$.
\end{remark}

In order to perform the rotation suggested in Remark~\ref{R:easy case},
we need two particular instances of 
Definition~\ref{D:maps}.
\begin{lemma} \label{L:fn}
For $N$ a nonnegative integer, let $f_N: A_L(0,1) \to A_L(0,1)$ be the map
fixing $L$ and pulling back $t$ to $t^{p^N}$.
Then for $\rho \in (\epsilon,1)$, 
$S(f_N^* \calE, \rho^{1/p^N}) \geq S(\calE,\rho)$.
Moreover, if $\del_i$ is dominant for $\calE_\rho$ for some
$i \neq n+1$, then $T(f_N^* \calE, \rho^{1/p^N}) = T(\calE,\rho)$.
\end{lemma}
\begin{proof}
The first assertion follows from Proposition~\ref{P:maps2}.
The second follows because if $\del_i$ is dominant for $\calE_\rho$
and $i \neq n+1$, then 
$T(f_N^* \calE, \rho^{1/p^N})$ and $T(\calE,\rho)$ can be computed
using the same formula.
\end{proof}

\begin{lemma} \label{L:gi}
Suppose $i \in \{1,\dots,n\}$ is such that
$\del_i$ is eventually dominant for $\calE$.
Let $g_i$ be the map given by Definition~\ref{D:maps} with
\[
u'_i = u_i + t, \qquad
u'_j = u_j \quad (j \neq i), \qquad
t' = t.
\]
Put $\calE' = g_i^* \calE$, and let $\calF'$ be the relativization of $\calE'$.
Let $b,b_{\rel}$ be the highest breaks of $\calE, \calF$. If $b>b_{\rel}+1$, then:
\begin{itemize}
\item
the highest break of $\calF'$ is $b-1$;
\item
for $\rho \in (0,1)$ sufficiently close to $1$,
the multiplicity of $\rho^{b-1}$ in $S(\calF',\rho)$
is the same as that of $\rho^b$ in $S(\calE,\rho)$.
\end{itemize}
\end{lemma}
\begin{proof}
The action of $\del_{n+1}$ on $g_i^* \calE$ is the pullback of the action
of $\del_{n+1} + \del_i$ on $\calE$, so the highest break of $\calF'$ is
the value of $b'$ satisfying 
\[
|\del_{n+1} + \del_i|_{\calE_\rho,\spect} = \rho^{-b'-1}
\]
for $\rho \in (0,1)$ sufficiently close to 1. For such $\rho$,
the spectral norms of $\del_i, \del_{n+1}$ on $\calE_\rho$ are $\rho^{-b},
\rho^{-b_{\rel}-1}$, respectively.
{}From this the claims are evident.
\end{proof}

\begin{lemma} \label{L:hard case}
Pick $i \in \{1, \dots, n+1\}$ such that $\del_i$ is eventually dominant
for $\calE$. Then at least one of the following statements is true.
\begin{itemize}
\item
For $\rho \in (0,1)$ sufficiently close to $1$, the scale multiset of
$\del_i$ on $\calE_\rho$ consists of a single element.
\item
There exists $\eta \in (0,1)$ such that $\calE$ is decomposable on
$A_L(\eta,1)$.
\end{itemize}
\end{lemma}
\begin{proof}
If $i = n+1$, then the claim follows by Remark~\ref{R:easy case},
so we assume $i \leq n$.
Let $b$ and $b_{\rel}$ be the highest breaks of $\calE$
and $\calF$,  respectively.
Assume that the first alternative does not hold; this forces $b > 0$.

Suppose to begin with that $b > b_{\rel} + 1$.
Put $\calE' = g_i^* \calE$ as in Lemma~\ref{L:gi}, 
and let $\calF'$ be the relativization of $\calE'$. 
Then $\calF'$ does not have a uniform break,
so by Lemma~\ref{L:cm}, it splits off a component of uniform break $b-1$.
We conclude that $\calE'$ is decomposable on some $A_L(\eta,1)$, as 
then is $\calE$, as desired.

In the general case, we can always pick
$N$ such that $b p^N > b_{\rel} + 1$. By Lemma~\ref{L:fn},
$f_N^* \calE$ has highest break $b p^N$, and
the first alternative of this lemma also does not hold for $f_N^* \calE$.
Moreover, by Proposition~\ref{P:frob scale}, the relativization of
$f_N^* \calE$ has highest break $b_{\rel}$.
We may thus apply the previous paragraph to split off a component of
$f_N^* \calE$ of highest break; since the splitting is unique, it descends
down the Galois group of the cover $f_N$ (after adjoining $p^N$-th roots of unity),
so $\calE$ is itself decomposable on some $A_L(\eta,1)$, as desired.
\end{proof}

\begin{proof}[Proof of Theorem~\ref{T:break decomp}]
It suffices to show that if $\calE$ is indecomposable over $A_L(\eta,1)$
for any $\eta \in (0,1)$ sufficiently close to 1, then $\calE$ has a uniform
break. This follows from Remark~\ref{R:easy case} if $\del_{n+1}$ is eventually 
dominant for $\calE$, and from Lemma~\ref{L:hard case} otherwise.
\end{proof}

It will be useful later to have a more uniform version of the rotation
construction used in Subsection~\ref{subsec:break decomp}, which comes
at the expense of enlarging the field $L$. 
(This generic rotation is inspired by the operation of \emph{generic
residual perfection} in \cite{borger}.)
The resulting construction will be used to study
the graded pieces of the ramification filtration.

\begin{prop} \label{P:uniform rotation}
Let $b$ be the highest break of $\calE$,
and suppose $b>1$.
Let $L'$ be the completion of $L(v_1, \dots, v_{n})$ for the
$(1,\dots,1)$-Gauss norm, viewed as a differential field of order $2n$ over
$K$. Let $\calE'$ be the pullback of $\calE$ along
the map $f: A_{L'}[0,1) \to A_L[0,1)$ given by
\[
f^*(u_i) = u_i^{p} + v_i t^{p-1} \quad (i=1,\dots,n),
\qquad
f^*(t) = \frac{t^p}{1-t^{p-1}}.
\]
Then $\calE'$ has highest break $p b - p+1$. In addition,
among the differentials 
\[
\frac{\del}{\del u_1}, \dots, \frac{\del}{\del u_n},
\frac{\del}{\del v_1}, \dots, \frac{\del}{\del v_{n}},
\frac{\del}{\del t},
\]
$\frac{\del}{\del t}$ (at least) is eventually dominant for $\calE'$.
\end{prop}
\begin{proof}
We first treat the case $n=0$. In this case, 
$g^*(t^{-1}) = t^{-p} - t^{-1}$, so this is an instance of
\cite[Lemma 5.13]{me-mono-over}.

In the general case, writing $\del'_1, \dots, \del'_{n+1}$ for the actions of
$\del_1, \dots, \del_{n+1}$ before the pullback, we have
\begin{align*}
\frac{\del}{\del u_i} &= pu_i^{p-1} \del'_{i} \\
\frac{\del}{\del v_i} &= t^{p-1} \del'_i \\
\frac{\del}{\del t} &= \frac{d}{dt} \left( \frac{t^p}{1-t^{p-1}} \right) \del'_{n+1} + 
\sum_{i=1}^n (p-1) v_i t^{p-2} \del'_i.
\end{align*}
We compute the scale of $\frac{\del}{\del t}$ by inspecting each
term separately: the contribution from $\del'_{n+1}$
can be treated as above, and the contribution from $\del'_i$ can
be treated directly
after invoking Proposition~\ref{P:maps3}.
This implies that the highest break of 
$\calE'$
is at least $pb-p+1$, with equality if and only if $\frac{\del}{\del t}$
is eventually dominant.

We compute the scale of $\frac{\del}{\del u_i}$ as if $u_i$ had pulled back
to $u_i^p$ and $t$ to $t^p$ (i.e., as for a Frobenius antecedent).
In particular, if $\frac{\del}{\del u_i}$ were eventually dominant
for $\calE'$, then the highest
break of $\calE'$ would be at most $b < pb-p+1$, contradiction.
Hence $\frac{\del}{\del u_i}$ is not eventually dominant.

We read off the scale of $\frac{\del}{\del v_i}$ directly:
it is eventually
dominant if and only if $\del'_i$ is, and in any case it cannot mask
$\frac{\del}{\del t}$. This proves the desired results.
\end{proof}

\begin{remark}
The calculations in this subsection may become more transparent when
checked against the examples produced by Artin-Schreier covers in
positive characteristic, as in Example~\ref{ex:artin-schreier}. Indeed,
many of these calculations were conceived with those examples
firmly in mind.
\end{remark}

\subsection{The differential Swan conductor}

Throughout this subsection, retain Hypothesis~\ref{H:solvable}.

\begin{defn}
By Theorem~\ref{T:break decomp}, there exists a multiset
$\{b_1, \dots, b_d\}$ such that
for all $\rho \in (0,1)$ sufficiently close to 1,
$S(\calE,\rho) = \{\rho^{b_1}, \dots, \rho^{b_d}\}$.
We call this multiset the 
\emph{break multiset} of $\calE$, denoted $b(\calE)$.
Define the \emph{(differential) Swan conductor} of $\calE$,
denoted $\Swan(\calE)$, as $b_1 + \cdots + b_d$.
\end{defn}

\begin{theorem} \label{T:Hasse-Arf}
The differential Swan conductor of $\calE$ is a nonnegative integer.
\end{theorem}
\begin{proof}
It suffices to check this in case $\calE$ is indecomposable over $A_L(\eta,1)$
for any $\eta \in (0,1)$ sufficiently close to 1.
Choose $i \in \{1, \dots, n+1\}$ such that $\del_i$ is eventually dominant for
$\calE$. By Lemma~\ref{L:hard case}, for $\rho \in (\epsilon,1)$ sufficiently
close to 1, the scale multiset of $\calE_\rho$
with respect to $\del_i$ consists of a single element. That means in the calculation
of the Newton polygon in Proposition~\ref{P:break piecewise}, the Newton polygon
must have only one slope, and so the integer $j$ 
can be taken to be $\rank(\calE)$. 
Consequently, the slopes of the function $f(r) = \log T(\calE,e^{-r})$
are always multiples of $1/\rank(\calE)$, as then is the highest break of
$\calE$. This proves the desired result.
\end{proof}

\begin{remark}
Proposition~\ref{P:maps3} implies that pulling $\calE$ along the map
$t \mapsto t^N$, for $N$ a positive integer coprime to $p$, has the effect
of multiplying $\Swan(\calE)$ by $N$. For Galois representations, this will
imply that the Swan conductor commutes appropriately with tamely
ramified base changes (Theorem~\ref{T:tame}).
\end{remark}

\begin{remark} \label{R:index}
In case $L=K$, one can interpret the integrality of $\Swan(\calE)$ by 
equating it to a certain local index \cite[Th\'eor\`eme~2.3-1]{cm4}. 
It would be interesting to
give a cohomological interpretation of our more general construction, perhaps
by relating it to an appropriate Euler characteristic.
\end{remark}

\begin{remark}
Liang Xiao points out that one can also prove Theorem~\ref{T:Hasse-Arf}
by reduction to the case of perfect residue field, for which one
may invoke Remark~\ref{R:index}. The argument is as follows. 
By Lemma~\ref{L:hard case},
we may assume that $\calE$ and its relativization
have respective uniform breaks $b, b_{\rel}$.
The perfect residue field case implies that $b_{\rel} \rank(\calE)$
is an integer. If $b \neq b_{\rel}$, we can
choose positive integers $m_1, m_2$ coprime to each other and to $p$
such that $m_i (b - b_{\rel}) > 1$ for $i=1,2$. 
If we pull back along $t \mapsto t^{m_i}$ and then apply
the rotation in Lemma~\ref{L:gi}, the highest break of the 
relativization becomes $m_i b - 1$, so $(m_i b - 1) \rank(\calE)$ is an
integer for $i=1,2$. This implies that $b \rank(\calE) \in \ZZ$.
\end{remark}

\section{Differential conductors for Galois representations}

In this section, we explain how to define differential Swan conductors for
certain $p$-adic Galois representations of complete discretely valued fields of
equal characteristic $p>0$ (including the discrete representations).
This uses a setup for turning representations into differential modules
due to Tsuzuki \cite{tsuzuki-amj}. For comments
on the mixed characteristic case, see Subsection~\ref{subsec:mixed}.

\subsection{Preliminaries: Cohen rings}

\begin{defn}
Let $k$ be a field of characteristic $p>0$. A \emph{Cohen ring} 
for $k$ is a complete discrete valuation ring $C_k$ with maximal ideal
generated by $p$, equipped with an isomorphism of its residue field
with $k$. \end{defn}

It can be shown that Cohen rings exist and are unique up to noncanonical
isomorphism; see \cite{bourbaki}. One can do better by carrying some extra 
data.

\begin{defn}
Define a \emph{based field} of characteristic $p>0$ to be a field
$k$ equipped with a distinguished $p$-basis $B_k$.
We view based fields as forming a category whose morphisms
from $(k,B_k)$ to $(k',B'_k)$ are morphisms $k \to k'$ of fields
carrying $B_k$ into $B'_k$.
\end{defn}

\begin{defn}
For $(k,B_k)$ a based field, a \emph{based Cohen ring} for $(k,B_k)$
is a pair $(C,B)$, where $C$ is a Cohen ring for $k$ and $B$ is a subset
of $C$ which lifts $B_k$.
\end{defn}

\begin{prop} \label{P:based Cohen}
There is a functor from based fields to based Cohen rings which is
a quasi-inverse of the residue field functor. In particular,
any map between based fields lifts uniquely to given based Cohen rings.
\end{prop}
\begin{proof}
This is implicit in Cohen's original paper \cite{cohen};
an explicit proof is given in \cite[Theorem~2.1]{whitney}.
Here is a sketch of another proof.
Let $W_n$ be the ring of $p$-typical Witt vectors of length $n$
over $k$, let $W$ be the inverse limit of the $W_n$,
let $F$ be the Frobenius on $W$,
and let $[\cdot]$ denote the Teichm\"uller map.
Put $B = \{[\overline{b}]: \overline{b} \in B_k\}$.
Let $C_n$ be the image of $F^n(W)[B]$ in $W_n$.
Then the projection $W_{n+1} \to W_n$ induces
a surjection $C_{n+1} \to C_n$.
Let $C$ be the inverse limit of the $C_n$; one then verifies that
$(C,B)$ is a based Cohen ring for $(k,B_k)$, and functoriality of the
construction follows
from functoriality of the Witt ring.
\end{proof}

\begin{remark} \label{R:whitney}
In fact, \cite[Theorem~2.1]{whitney} asserts something slightly stronger:
if $(C,B)$ is a based Cohen ring of $(k,B_k)$, $R$ is any complete local ring 
with residue field $k$, and $B_R$ is a lift of $B_k$ to $R$, then
there is a unique ring homomorphism $C \to R$ inducing the identity on $k$
and carrying $B$ to $B_R$.
\end{remark}

\subsection{Galois representations and $(\phi, \nabla)$-modules}

\begin{hypothesis}
For the remainder of this section,
let $R$ be a complete discrete valuation ring of
equal characteristic $p>0$, with fraction field $E$ and residue field $k$.
Let $k_0 = \cap_{n \geq 0} k^{p^n}$ be the maximal perfect subfield
of $k$;  note that $k_0$ embeds canonically into $R$ (whereas if $k \neq k_0$,
then $k$ embeds but not canonically).
\end{hypothesis}

\begin{convention}
Put $G_E = \Gal(E^{\sep}/E)$.
Let $\calO$ be the integral closure of $\ZZ_p$ in a finite extension
of $\QQ_p$, whose residue field $\FF_q$ is contained in $k$.
Throughout this section, a ``representation'' will be a 
continuous representation $\rho: G_E \to \GL(V)$, where $V = V(\rho)$
is a finite free $\calO$-module. (One can also consider representations
on finite dimensional $\Frac(\calO)$-vector spaces, by choosing lattices;
for brevity, we stick to statements for integral representations,
except for Remark~\ref{R:change lattice}.)
\end{convention}

\begin{defn}
Fix a based Cohen ring $(C_E, B)$ with residue field $E$;
note that $C_E$ is canonically a
$W(\FF_q)$-algebra. Put 
\[
\Gamma = C_E \otimes_{W(\FF_q)} \calO.
\]
Let $\Omega^1_{\Gamma/\calO}$ be the completed (for the $p$-adic topology)
direct sum of $\Gamma\,db$ over all
$b \in B$, i.e., the inverse limit over $n$ 
of $\oplus_{b \in B} (\Gamma/p^n \Gamma)\,db$;
then there is a canonical derivation
$d: \Gamma \to \Omega^1_{\Gamma/\calO}$. Note that all of this data
stays canonically independent of the choice of $B$ as long as 
$C_E$ remains fixed.
\end{defn}

\begin{defn}
A \emph{$\nabla$-module} over $\Gamma$ is a finite free $\Gamma$-module $M$
equipped with an integrable connection $\nabla: M \to M \otimes_\Gamma
\Omega^1_{\Gamma/\calO}$; integrability means that the composition of
$\nabla$ with the map $M \otimes \Omega^1_{\Gamma/\calO} \to
M \otimes \wedge^2_\Gamma \Omega^1_{\Gamma/\calO}$ induced by $\nabla$
is the zero map.
\end{defn}

\begin{defn}
A \emph{Frobenius lift} on $\Gamma$ is an endomorphism 
$\phi: \Gamma \to \Gamma$ fixing $\calO$ and lifting the $q$-power
Frobenius map on $E$. For instance, there is a unique such $\phi$
carrying $b$ to $b^q$ for each $b \in B$ (induced by the Frobenius action
on the construction given in Proposition~\ref{P:based Cohen}); 
we call this $\phi$ the 
\emph{standard Frobenius lift with respect to $B$}.
A \emph{$\phi$-module} (resp.\ \emph{$(\phi, \nabla)$-module}) over $\Gamma$ 
is a finite free module (resp.\ $\nabla$-module) $M$ over $\Gamma$
equipped with an isomorphism $F: \phi^* M \cong M$ of
modules (resp.\ of $\nabla$-modules); we interpret $F$ as a semilinear
action of $\phi$ on $M$.
\end{defn}

\begin{defn}
For any representation $\rho$, put
\[
D(\rho) = (V(\rho) \otimes_{\calO} \widehat{\Gamma^{\unr}})^{G_E}.
\]
By Hilbert's Theorem 90, the natural map
\[
D(\rho) \otimes_{\Gamma} \widehat{\Gamma^{\unr}} \to V(\rho) \otimes_{\calO}
\widehat{\Gamma^{\unr}}
\]
is a bijection; in particular, $D(\rho)$ is a free $\Gamma$-module
and $\rank_\Gamma(D(\rho))=\rank_{\calO}(V)$.
If we equip $\Gamma^{\unr}$ and its completion
with actions of the derivation $d$ and any Frobenius lift $\phi$
(acting  trivially on $V(\rho)$),
we obtain by restriction a Frobenius action and connection on $D(\rho)$,
turning it into a $(\phi, \nabla)$-module.
\end{defn}

\begin{prop} \label{P:d equivalence}
For any Frobenius lift $\phi$ on $\Gamma$,
the functor $D$ from representations to $\phi$-modules over $\Gamma$ is
an equivalence of categories.
\end{prop}
\begin{proof}
Given a $\phi$-module $M$ over $\Gamma$, put
\[
V(M) = (M \otimes_\Gamma \widehat{\Gamma^{\unr}})^{\phi = 1}.
\]
As in \cite[A1.2.6]{fontaine} or \cite[Theorem~4.1.3]{tsuzuki-amj},
one checks that $V$ is a quasi-inverse to $D$.
\end{proof}

\begin{prop}
For any Frobenius lift $\phi$ on $\Gamma$,
any $\phi$-module over
$\Gamma$ admits a unique structure of $(\phi, \nabla)$-module.
Consequently, 
the functor $D$ from representations to $(\phi, \nabla)$-modules 
over $\Gamma$ is
an equivalence of categories.
\end{prop}
\begin{proof}
Existence of such a structure follows from Proposition~\ref{P:d equivalence},
so we focus on uniqueness. Let $M$ be a $(\phi, \nabla)$-module over $\Gamma$.
Let $\frac{\del}{\del b}$ be the derivations dual to the $db$ for $b \in B$.
Let $\be_1, \dots, \be_m$ be a basis of $M$, and let $\Phi$ and $N_b$
be the matrices via which $\phi$ and $\frac{\del}{\del b}$ act on this basis.
Then the fact that the $\phi$-action on $M$ respects the $\nabla$-module 
structure implies that
\begin{equation} \label{eq:determine nabla}
N \Phi + \frac{\del \Phi}{\del b} = \frac{\del \phi(b)}{\del b} \Phi N.
\end{equation}
Let $\pi$ be a uniformizer of $\calO$; note that
$\frac{\del \phi(b)}{\phi b} \equiv 0 \pmod{\pi}$ because
$\phi(b) \equiv b^q \pmod{\pi}$. Consequently, for fixed $\Phi$,
if $N_b$ is uniquely
determined modulo $\pi^m$, then the right side of
\eqref{eq:determine nabla} is determined modulo $\pi^{m+1}$, as then
is $N_b \Phi$. Since $\Phi$ is invertible, $N_b$ is also
determined modulo $\pi^{m+1}$. By induction, $N_b$ is uniquely determined
by $\Phi$ for each $b$, as desired.
\end{proof}

\subsection{Representations with finite local monodromy}

We now distinguish the class of representations for which we define
differential Swan conductors.

\begin{defn}
Let $I_E = \Gal(E^{\sep}/E^{\unr})$ be the inertia subgroup of $G_E$.
We say a representation $\rho$ \emph{has finite local monodromy} if the 
image of $I_E$ under $\rho$ is finite.
\end{defn}

For representations with finite local monodromy, we can refine the construction
of the $(\phi, \nabla)$-module associated to $\rho$.

\begin{hypothesis} \label{H:lifted cohen}
For the remainder of this subsection,
assume that $k$ admits a \emph{finite} $p$-basis.
Assume also that the based Cohen ring $(C_E, B)$ has been chosen with $B = 
B_0 \cup \{t\}$, where $t$ lifts a uniformizer of $E$, and $B_0$ lifts
elements of $R$ which in turn lift a $p$-basis of $k$.
\end{hypothesis}

\begin{defn}
By the proof of the Cohen structure theorem, or by Remark~\ref{R:whitney},
there is 
a unique embedding of $k$ into $R$ whose image contains the image
of $B_0$ under reduction to $E$.
Applying Proposition~\ref{P:based Cohen} to the map $k \to R$,
we obtain an embedding of a Cohen ring $C_k$ for $k$ into $C_E$, the
image of which contains $B_0$.
Put $\calO_k = C_k \otimes_{W(\FF_q)} \calO$. Then 
each $x \in \Gamma$ can be written formally as a 
sum $\sum_{i \in \ZZ} x_i t^i$ with $x_i \in \calO_k$, such that for each
$n$, the indices $i$ for which $v_{\calO_k}(x_i) \leq n$
are bounded below.
For $n$ a nonnegative integer, we define the \emph{partial valuation function}
$v_n: \Gamma \to \ZZ \cup \{\infty\}$ by
\[
v_n(x) = \min\{i \in \ZZ: v_{\calO_k}(x_i) \leq n\}.
\]
For $r>0$, put
\[
\Gamma^r = \{x \in \Gamma: \lim_{n \to \infty} v_n(x) + rn = \infty\};
\]
this is a subring of $\Gamma$. Put $\Gamma^\dagger = \cup_{r>0} \Gamma^r$;
we may speak of $\nabla$-modules over $\Gamma^\dagger$ using the
same definition as for $\Gamma$, using for the module of differentials
\[
\Omega^1_{\Gamma^\dagger/\calO_k}
= \bigoplus_{b \in B} \Gamma^\dagger\,db.
\]
(Here we are using the finiteness of the
$p$-basis to avoid having to worry about a completion.)
 If $\phi$ is a Frobenius lift
carrying $\Gamma^\dagger$ into itself, we may also
define $\phi$-modules and $(\phi, \nabla)$-modules
over $\Gamma^\dagger$ as before.
\end{defn}

\begin{defn}
Since $\calO_k \subset \Gamma^\dagger$, we can identify a copy of
$\calO_k^{\unr}$ inside $(\Gamma^{\dagger})^{\unr}$. Using this identification,
put
\[
\tilde{\Gamma}^\dagger = \widehat{\calO_k^{\unr}} \otimes_{\calO_k^{\unr}}
(\Gamma^{\dagger})^{\unr} \subset \widehat{\Gamma^{\unr}}.
\]
For $\rho$ a representation, put
\begin{align*}
D^\dagger(\rho) &= D(\rho) \cap (V(\rho) 
\otimes_\calO \tilde{\Gamma}^\dagger) \\
&= (V(\rho) \otimes_\calO \tilde{\Gamma}^\dagger)^{G_E}.
\end{align*}
Again, $D^\dagger(\rho)$ inherits a connection, and an action of any
Frobenius lift $\phi$ acting on $\Gamma^\dagger$.
Note that the natural map
\[
(D^\dagger(\rho) \otimes_{\Gamma^\dagger} \tilde{\Gamma}^\dagger)
\to (V(\rho) \otimes_{\calO} \tilde{\Gamma}^\dagger)
\]
is always injective, and it is surjective if and only if $\rho$
has finite local monodromy.
\end{defn}

The following is essentially \cite[Theorem~3.1.6]{tsuzuki-amj}.
\begin{prop} \label{P:full faith}
Let $\phi$ be a Frobenius lift on $\Gamma$ acting on $\Gamma^\dagger$.
The base change functor from $(\phi, \nabla)$-modules over $\Gamma^\dagger$
to $(\phi, \nabla)$-modules over $\Gamma$ is fully faithful.
\end{prop}
\begin{proof}
Using internal Homs, we may rephrase this as follows: if $M$ is
a $(\phi, \nabla)$-module over $\Gamma^\dagger$, then
\[
(M \otimes \Gamma)^{\phi=1,\nabla=0} \subset M.
\]
In particular, it is sufficient to check this using only the $dt$ component
of the connection. In this case, we may replace $\Gamma$ by the completion
of $\Gamma[\phi^{-n}(b): b \in B_0, n \in \ZZ_{\geq 0}]$, to get into the
case where $R$ has perfect residue field. We may then conclude by applying
a result of Tsuzuki \cite[4.1.3]{tsuzuki-etale}.
\end{proof}

The following is essentially \cite[Theorem~4.2.6]{tsuzuki-amj}.
\begin{theorem} \label{T:dagger equiv}
Let $\phi$ be a Frobenius lift on $\Gamma$ acting on $\Gamma^\dagger$.
Then $D^\dagger$ and restriction induce
 equivalences between the following categories:
\begin{enumerate}
\item[(a)] representations with finite
local monodromy;
\item[(b)]
$(\phi, \nabla)$-modules over $\Gamma^\dagger$;
\item[(c)]
$\nabla$-modules over $\Gamma^\dagger$ equipped with $\phi$-actions
over $\Gamma$.
\end{enumerate}
In particular, if a $\nabla$-module over $\Gamma^\dagger$ admits
a $\phi$-action over $\Gamma$, that action is defined already over
$\Gamma^\dagger$.
\end{theorem}
\begin{proof}
The functor from (a) to (b) is $D^\dagger$, while the functor from (b) to
(c) is restriction. The functor from (c) back to (a) will be induced by
$V$; once it is shown to be well-defined, it will be clear
that the three functors compose to the identity
starting from any point.

To obtain the functor from (c) to (a), we must prove 
that if $M$ is a $\nabla$-module over $\Gamma^\dagger$ such that
$M \otimes \Gamma$ admits a compatible $\phi$-action, 
then the corresponding representation $V(M)$ has finite local monodromy.
It suffices to check this after replacing $E$ by a finite extension,
which can be chosen to ensure the existence of an isomorphism
$(M/2pM) \otimes \Gamma \cong (\Gamma/2p\Gamma)^m$ of
$\phi$-modules. In this case we claim that $V(M)$ is actually
unramified; 
as in the proof of Proposition~\ref{P:full faith},
it suffices to check this using only the $dt$ component of $\nabla$,
and hence to reduce to the case of $R$ having perfect residue field.
This case is treated by the proof of \cite[Proposition~5.2.1]{tsuzuki-amj},
but not by its statement (which requires a $\phi$-action over 
$\Gamma^\dagger$); for a literal citation, see
\cite[Proposition~4.5.1]{me-fake}.
\end{proof}

\subsection{$(\phi,\nabla)$-modules over $\calR$}
\label{subsec:phi-nabla}

Throughout this subsection, retain
Hypothesis~\ref{H:lifted cohen},
and write $L$ for $\Frac(\calO_k)$ and
$\calR$ for $\calR_L$. 
The choices made so far determine
an embedding $\Gamma^\dagger \hookrightarrow \calR_L$, and any
Frobenius $\phi$ acting on $\Gamma^\dagger$ extends continuously
to $\calR_L$ (as in \cite[\S 2]{me-local}).
We may thus define $\phi$-modules, $\nabla$-modules,
and $(\phi,\nabla)$-modules over $\calR$
using the same definitions as over $\Gamma$. 

\begin{remark} \label{R:geometrize}
{}From a $\nabla$-module over $\calR$, we may construct a 
$\nabla$-module on $A_L(\epsilon,1)$ for some
$\epsilon \in (0,1)$. The construction is unique
in the following sense: any two such $\nabla$-modules become
isomorphic on $A_L(\eta,1)$ for some $\eta \in (0,1)$.
Conversely, since any locally free sheaf on $A_L(\eta,1)$ is freely
generated by global sections (because $L$ is spherically complete;
see for instance \cite[Theorem~3.14]{me-mono-over}),
any $\nabla$-module on $A_L(\eta,1)$ gives rise to a 
$\nabla$-module over $\calR$.
\end{remark}

Remark~\ref{R:geometrize} is sufficient for the construction
of the differential Swan conductor associated to a representation
of finite local monodromy. However, for completeness, we
record some related facts, including the analogue of the $p$-adic
local monodromy theorem.

\begin{lemma} \label{L:descend nabla}
Let $M$ be a $\phi$-module over $\Gamma^\dagger$ such that
$M \otimes \calR$ admits the structure of a $(\phi, \nabla)$-module.
Then this structure is induced by a $(\phi, \nabla)$-module structure
on $M$ itself, and so
$M$ corresponds to a representation with finite local monodromy.
\end{lemma}
\begin{proof}
By \cite[Proposition~7.1.7]{me-slope}, the action of $\frac{\del}{\del t}$
on $M \otimes \calR$ acts on $M$ itself. 
Also, for any $b \in B_0$, we may change the $p$-basis by replacing $b$ by
$b + t$, and then the same argument shows that the action of 
$\frac{\del}{\del t} + \frac{\del}{\del b}$ on $M \otimes \calR$ acts
on $M$ itself. (This is another instance of rotation in the sense of
Remark~\ref{R:easy case}.) 
We conclude that $\nabla$ itself acts on $M$, so
we may apply Theorem~\ref{T:dagger equiv} to conclude.
\end{proof}

\begin{defn}
A $\phi$-module (resp.\ $(\phi,\nabla)$-module)
$M$ over $\calR$ is \emph{unit-root} if 
it has the form $M_0 \otimes \calR$ for some
$\phi$-module (resp.\ $(\phi, \nabla)$-module) $M_0$ over
$\Gamma^\dagger$.
By Lemma~\ref{L:descend nabla}, a $(\phi,\nabla)$-module over $\calR$
is unit-root if and only if its underlying $\phi$-module is unit-root.
\end{defn}

\begin{prop}
The base extension functor from the isogeny category of 
unit-root $\phi$-modules over $\Gamma^\dagger$ (i.e.,
$\phi$-modules over $\Gamma^{\dagger}[\frac 1p]$ obtained by
base extension from $\Gamma^\dagger$) to 
the category of unit-root $\phi$-modules over $\calR$
is an equivalence of categories.
\end{prop}
\begin{proof}
This is \cite[Theorem~6.2.3]{me-slope}.
\end{proof}

\begin{defn}
Let $s = c/d$ be a rational number written in lowest terms.
A $\phi$-module (resp.\ $(\phi, \nabla)$-module) $M$ over $\calR$ is
\emph{pure} (or \emph{isoclinic}) of slope $s$ if
there exists a scalar $\lambda \in K^*$ of valuation $c$
such that the $\phi^d$-module (resp.\ $(\phi^d, \nabla)$-module)
obtained from $M$ by twisting the $\phi^d$-action by $\lambda^{-1}$
is 
unit-root. In particular, by Theorem~\ref{T:dagger equiv}, the
$\nabla$-module structure on $M$ corresponds to a representation
with finite local monodromy after replacing $\calO$ by a finite extension.
\end{defn}

\begin{theorem} \label{T:slope filt}
Let $M$ be a $\phi$-module (resp.\ $(\phi, \nabla)$-module) $M$ over $\calR$.
Then there exists a unique filtration $0 = M_0 \subset M_1 \subset \cdots
\subset M_l = M$ of $M$ by saturated $\phi$-submodules (resp.\
$(\phi,\nabla)$-submodules) such that each quotient $M_i/M_{i-1}$ is pure
of some slope $s_i$ as a $\phi$-module, and $s_1 < \cdots < s_l$.
\end{theorem}
\begin{proof}
In the $\phi$-module case, this is \cite[Theorem~6.10]{me-local}
or \cite[Theorem~6.4.1]{me-slope}. In the $(\phi, \nabla)$-module case,
it suffices to check that the filtration of the underlying $\phi$-module
is respected by $\nabla$. For this, we proceed
as in \cite[Theorem~7.1.6]{me-slope}: 
for each derivation $\frac{\del}{\del b}$,
we get a
morphism of $\phi$-modules $M_1 \to (M/M_1) \otimes \calR\,db$.
The former is pure of slope $s_1$, whereas the latter admits a slope filtration
in which each slope is strictly greater than $s_1$ (the slope of
$\calR\,db$ being positive). By \cite[Proposition~4.6.4]{me-slope}, 
that morphism is zero, proving that $M_1$ is respected by each derivation.
Hence $M_1$ is a $(\phi,\nabla)$-submodule, and repeating the argument
on $M/M_1$ yields the claim.
\end{proof}

\begin{remark}
One may apply Theorem~\ref{T:dagger equiv} to each individual quotient
of the filtration produced by Theorem~\ref{T:slope filt}.
(Alternatively, one may project $\nabla$ onto the $dt$ component and
directly invoke the $p$-adic local monodromy theorem; this allows
the invocation of \cite{andre} or \cite{mebkhout} in place of
\cite{me-local}.) It is an interesting question, which we have not
considered, whether one can show that the category of
$(\phi, \nabla)$-module $M$ over $\calR_L$ is equivalent to a category
of representations of $G_E$ times an algebraic group 
over $\Frac(\calO)$, as in \cite[Theorem~4.45]{me-mono-over}.
\end{remark}

\subsection{Defining the differential Swan conductor}

In order to use Theorem~\ref{T:dagger equiv} to define the differential
Swan conductor of a representation $\rho: G_E \to \GL(V)$ with finite
local monodromy, we must check that the answer 
does not depend on the auxiliary choices
we made along the way. 
(Note that the choice of $\phi$ does not matter: it is
only used to define the Frobenius action on $D^\dagger(\rho)$, whereas only
the connection is used to compute the conductor.)

\begin{prop} \label{P:indep of choices}
Suppose that $k$ admits a finite $p$-basis.
For $\rho$ a representation with finite local monodromy,
the isomorphism type of the $\nabla$-module $D^\dagger(\rho)$
does not depend on the choice of the Cohen ring $C_E$ or the
lifted $p$-basis $B$. 
\end{prop}
\begin{proof}
By Proposition~\ref{P:based Cohen}, the construction of $C_E$ is functorial
in pairs $(E, \overline{B})$, where $\overline{B}$ is a $p$-basis of
$E$. It thus suffices to check that if for $i=1,2$,
$\overline{B_i}$ is a $p$-basis of $E$ consisting of a uniformizer 
$\overline{t_i}$ of $R$ and a lift $\overline{B_{i,0}}$ to $R$
of a $p$-basis of $k$ over $k_0$, then the modules $D^\dagger(\rho)$
constructed using lifts of $\overline{B_1}$ and $\overline{B_2}$
are isomorphic, compatibly with some isomorphism of the underlying
rings $\Gamma^\dagger$.

Let $(C_E,B_1)$ be a based Cohen ring lifting $(E,\overline{B_1})$;
write $C_{k,1}, t_1$ instead of $C_k,t$.
Define $B_{2,0}$ by choosing, for each $\overline{b} \in B_{2,0}$, 
a lift $b$ of $\overline{b}$ in $C_{k,1} \llbracket t_1 \rrbracket$. 
Then choose
$t_2$ to be a lift of $\overline{t_2}$ belonging to
$t_1 C_{k,1} \llbracket t_1 \rrbracket$. We can then view
$(C_E, B_2)$ as a based Cohen ring lifting $(E, \overline{B_2})$,
containing a Cohen ring $C_{k,2}$ for $k$.

Since we used the same ring $C_E$ for both lifts, we may identify
the two rings $\Gamma$. Although $C_{k,1} \neq C_{k,2}$ in general,
we did ensure by construction that
$C_{k,1} \llbracket t_1 \rrbracket = C_{k,2} \llbracket t_2 \rrbracket$.
Consequently, the two rings
$\Gamma^\dagger$ constructed inside $\Gamma$ coincide,
and we may identify the two copies of $\tilde{\Gamma}^\dagger$.
This gives an identification of the two modules $D^\dagger(\rho)$,
as desired.
\end{proof}

\begin{defn} \label{D:rep diff}
Suppose to start
that $k$ is finite over $k^p$, i.e., any $p$-basis of $k$
or of $E$ is finite.
For $\rho: G_E \to \GL(V)$ a representation with finite local monodromy,
with $V$ a finite dimensional $\calO$-module, we may now define 
the \emph{differential highest break},
\emph{differential break multiset}, and 
\emph{differential Swan conductor} by constructing the 
$(\phi, \nabla)$-module $D^\dagger(\rho)$, for some Cohen ring $C_E$
and some lifted $p$-basis $B$, and computing the corresponding
quantities associated to the underlying $\nabla$-module of $D^\dagger(\rho)$
tensored with the Robba ring $\calR_{\Frac(\calO_k)}$ (as in 
Remark~\ref{R:geometrize}). By
Proposition~\ref{P:maps1} (to change the $p$-basis of $k$)
and Proposition~\ref{P:indep of choices}, this definition depends only
on $\rho$ and not on any auxiliary choices.
For general $k$, we may choose a finite subset $B_1$ of $B$ containing a lift
$t$ of a uniformizer of $R$, project onto the span of the $db$ for $b \in B_1$,
and compute a conductor that way; 
this has the same effect as passing from $E$ to $E_1 = \widehat{E_0}$, where
\[
E_0 = E(\overline{b}^{1/p^n}: b \in B \setminus B_1,
n \in \NN_0).
\]
We define the differential Swan conductor of $\rho$ in
this case to be the supremum over all choices of $B$ and $B_1$;
it will turn out to be finite (Corollary~\ref{cor:compare2})
and hence integral by Theorem~\ref{T:Hasse-Arf}.
\end{defn}

\begin{defn} \label{D:rectifying}
Let $E'/E$ be a finite separable extension, let $\overline{B}$
be a $p$-basis of $E$ containing a uniformizer $\overline{t}$ of $E$,
and put $\overline{B_0} = \overline{B} \setminus \{\overline{t}\}$. We say a
subset $\overline{B_2}$ of $\overline{B_0}$ 
is a \emph{rectifying set} for $E'/E$ if,
putting $E_2 = \widehat{E_0}$ for
\[
E_0 = E(\overline{b}^{1/p^n}: \overline{b} \in \overline{B_2},
n \in \NN_0),
\]
the extension $(E' \otimes_E E_2) / E_2$ has separable residue field 
extension. Beware that it is not enough for the residue field of
$E_2$ to contain the perfect closure of $k$
in the residue field of $E'$. For instance, if $p>2$, $\overline{b_1}, 
\overline{b_2} \in
\overline{B}$, and 
\[
E' = E[z]/(z^p - z - \overline{b_1} \overline{t}^{-2p} 
- \overline{b_2} \overline{t}^{-p}),
\]
then $E'$ has residue field $k(\overline{b_1}^{1/p})$,
but $\overline{B_2} = \{\overline{b_1}\}$ is not a rectifying set
because the residue field of $E' \otimes_E E_2$ contains 
$\overline{b_2}^{1/p}$.
\end{defn}

\begin{lemma} \label{L:finite rectify}
With notation as in Definition~\ref{D:rectifying},
$\overline{B}$ contains a finite rectifying set for $E'/E$.
\end{lemma}
\begin{proof}
Use $\overline{B}$ to embed $k$ into $E$.
By induction on the degree of $E'/E$, we may reduce to the case of 
an Artin-Schreier extension
\[
E' = E[z]/(z^p - z - a_n \overline{t}^{-n} - \cdots - a_1 
\overline{t}^{-1} - a_0)
\]
with $a_i \in k$. In this case, pick any $N \in \NN_0$ with $p^N > n$,
and write each $a_i$ as a $k^{p^N}$-linear
combination of products of powers of elements of $\overline{B_0}$.
Only finitely many elements of 
$\overline{B_0}$ get used; those form a rectifying set.
\end{proof}

\begin{prop} \label{P:reconcile}
Suppose that there exists a finite separable extension $E'$ of $E$
whose residue field is separable over $k$, such that $\rho$ is unramified
on $G_{E'}$. Then the differential break multiset and
Swan conductor of a representation
$\rho$ with finite local monodromy can be computed
with respect to $\{t\}$, and it agrees with the usual break
multiset and Swan conductor.
\end{prop}
\begin{proof}
It suffices to consider $\rho$ irreducible and check equality for the
highest breaks.
Note that the usual highest break is insensitive to further residue
field extension, because it can be computed using Herbrand's formalism
as in \cite[Chapter~IV]{serre}. It thus agrees with the
differential highest break computed with respect to $\{t\}$: namely,
this claim reduces to the case where $k$ is perfect,
for which  see \cite[Theorem~5.23]{me-mono-over} and references thereafter.

It remains to show that for any $B_1$, $\frac{\del}{\del t}$ must be
eventually dominant. Suppose the contrary, and pick $b \in B_1 
\setminus \{t\}$ such that $\frac{\del}{\del b}$ is eventually dominant.
By a tame base change (invoking Proposition~\ref{P:maps3}), 
we can force the gap between the differential
highest breaks computed with respect to $B_1$ and with respect to
$\{t\}$ to be greater 
than 1; then a rotation as in Lemma~\ref{L:hard case} sending $b$ to $b+t$
raises the differential highest break
computed with respect to $\{t\}$. But that contradicts the previous paragraph:
both before and after rotation, the differential highest break
computed with respect to $\{t\}$ must coincide with the usual highest break.

We deduce that $\frac{\del}{\del t}$ is eventually dominant, proving the
claim.
\end{proof}
\begin{cor} \label{cor:compare1}
In the notation of Definition~\ref{D:rep diff},
suppose that there exists a finite separable extension $E'$ of $E$
such that $\rho$ is unramified on $G_{E'}$, and that
the image of $B_1 \setminus \{t\}$ in $E$ is a rectifying set for $E'/E$. Then
the differential Swan conductor of $\rho$ computed using 
$(B \setminus B_1) \cup \{t\}$ is equal to that computed using $t$.
\end{cor}
\begin{cor} \label{cor:compare2}
In the notation of Definition~\ref{D:rep diff},
suppose that there exists a finite separable extension $E'$ of $E$
such that $\rho$ is unramified on $G_{E'}$, and that
the image of $B_1 \setminus \{t\}$ in $E$ is a rectifying set for $E'/E$. Then
the differential Swan conductor of $\rho$ 
is equal to that computed using $B_1$.
\end{cor}

For completeness, we record the following observations.
\begin{theorem}
The differential Swan conductor of any representation with finite
local monodromy is a nonnegative integer.
\end{theorem}
\begin{proof}
By Lemma~\ref{L:finite rectify} and
Corollary~\ref{cor:compare2}, the conductor can be computed using a
finite set $B_1$; we may thus apply Theorem~\ref{T:Hasse-Arf}.
\end{proof}

\begin{theorem} \label{T:tame}
Let $E'$ be a tamely ramified extension of $E$ of ramification degree $m$.
Let $\rho$ be a representation of $G_E$ with finite local monodromy,
and let $\rho'$ be the restriction of $\rho$ to $G_{E'}$. Then
$\Swan(\rho') = m \Swan(\rho)$.
\end{theorem}
\begin{proof}
Apply Proposition~\ref{P:maps3}.
\end{proof}

\begin{example} \label{ex:artin-schreier}
As an example, consider a nontrivial character of the Artin-Schreier
extension $E[z]/(z^p - z - \overline{x})$. 
The corresponding differential module will be a Dwork isocrystal, i.e.,
a rank one $\nabla$-module with generator $\bv$ such that
\[
\nabla(\bv) = \pi \bv \otimes dx,
\]
for $x$ some lift of $\overline{x}$ and $\pi$ a $(p-1)$-st root of 
$-p$.
One computes that the differential
Swan conductor for this character is equal to the least integer
$m \geq 0$ such that $v_E(x - y^p + y) \geq -m$ for all $y \in E$.
This agrees with the definition given by Kato \cite{kato} of the
Swan conductor of a character; note that the
conductor is allowed to be divisible by $p$ if and only if $k$ is imperfect.
\end{example}

\begin{remark} \label{R:change lattice}
Given a representation $\rho: G_E \to \GL(V)$, where $V$ is a finite
dimensional $\Frac(\calO)$-vector space, we may define a differential
Swan conductor for it by picking a $\rho$-stable $\calO$-lattice
of $V$ and proceeding as in Definition~\ref{D:rep diff}.
Changing the lattice will not change the resulting $\nabla$-module 
over $\calR_{\Frac(\calO_k)}$, so we get a well-defined numerical invariant
of $\rho$ also.
\end{remark}

Defining conductors for Galois representations is tantamount to filtering
the Galois group; let us now make this explicit.
\begin{defn}
Put $G_E^0 = I_E$. For $r> 0$, let $R_r$ be the set of representations $\rho$
with highest break less than $r$, and put
\[
G_E^r = \bigcap_{\rho \in R_r} (I_E \cap \ker(\rho)).
\]
Note that $\rho \in R_r$ if and only if $G_E^r \subseteq I_E \cap \ker(\rho)$;
this reduces to the fact that $R_r$ is stable under tensor product
and formation of subquotients.
We call $G_E^r$ the \emph{differential upper numbering filtration} on $G_E$.
Write $G_E^{r+}$ for the closure of $\cup_{s>r} G_E^s$; note that
$G_E^r = G_E^{r+}$ for $r$ irrational, because differential highest breaks
are always rational numbers.
\end{defn}

As in the perfect residue field case, the graded pieces of the
upper numbering filtration are particularly simple.
\begin{theorem} \label{T:upper numbering}
For $r > 0$ rational, $G_E^r/G_E^{r+}$ is abelian and killed by $p$.
\end{theorem}
\begin{proof}
Let $E'$ be a finite Galois extension of $E$ with $\Gal(E'/E) = G$;
then we obtain an induced filtration on $G$ by taking
$G^r$ to be the image of $G_E^r$ under the surjection $G_E \to G$.
It suffices to check that $G^r/G^{r+}$ is abelian and killed by $p$;
moreover, we may quotient further to reduce to the case where $G^{r+}$
is the trivial group but $G^r$ is not. 
Let $\rho$ be the regular representation of $G$; then $\rho$ has highest
break $r$. Let $S$ be the set of irreducible constituents of
$\rho$ of highest break strictly less than $r$; we are then trying to show
that the intersection of $\ker(\psi) \subseteq G$ over all $\psi \in S$
is an elementary abelian $p$-group.

By Corollary~\ref{cor:compare2}, we may reduce to the case where
the lifted $p$-basis $B$ of Hypothesis~\ref{H:lifted cohen} is finite;
put $B_0 = \{b_1, \dots, b_n\}$. 
By making a tame base change, we can force all nonzero ramification
breaks to be greater than 1. By another base change
(passing from $\Frac C_E$ to the completion of
$C_E(v_1, \dots, v_{n})$ for the $(1, \dots, 1)$-Gauss norm), 
we can add extra elements $v_1, \dots, v_{n}$ to $B$, then
perform the operation described in
Proposition~\ref{P:uniform rotation}. 
Each nonzero ramification break $m$ before
the operation corresponds to the break $p 
m - p + 1$ afterwards, so the desired result
may be checked afterwards. But now $\frac{\del}{\del t}$ is dominant
on every irreducible component of $\rho$,
so we may reduce to the case of perfect residue field
and (by Proposition~\ref{P:reconcile}) the usual upper numbering filtration.
In this case, the claim is standard: it
follows from the fact that the upper numbering
filtration can be constructed by renumbering the lower numbering filtration
\cite[\S IV.3]{serre}, for which the claim is easy to check
\cite[\S IV.2, Corollary~3 of Proposition~7]{serre}.
\end{proof}

\begin{remark}
Note that the definition of the differential Swan conductor of a
representation is invariant
under enlarging $\calO$, because the differential Swan conductor of
a $\nabla$-module is invariant under enlarging the constant field $K$.
\end{remark}

\subsection{Reconciliation questions}

By introducing a numerical invariant of representations and calling it
a conductor, one begs various reconciliation questions with other definitions.
To begin with, it is known (and was a motivation of our construction) that
in the traditional case of a perfect residue field, one computes the
right numbers; see Proposition~\ref{P:reconcile}.

In the general case, there is a definition of the ``logarithmic
conductor'' due to Abbes and Saito \cite{abbes-saito1, abbes-saito2}.
Following Matsuda \cite{matsuda-dwork}, one is led to ask the following.
\begin{question} \label{Q:abbes-saito}
For $\rho$ a representation with finite local monodromy, does the
differential Swan conductor agree with the Abbes-Saito logarithmic
conductor in equal characteristic?
\end{question}
It is easy to check the affirmative answer for Artin-Schreier characters.
An affirmative answer in the general case would have the beneficial consequences
of verifying the Hasse-Arf theorem for the Abbes-Saito
conductor in equal characteristic. Some progress on this question has
been made recently by Bruno Chiarellotto and Andrea Pulita, and independently
by Liang Xiao.

One might also try to reconcile our definition with conductors
for Galois representations over a two-dimensional local field, as in
Zhukov \cite{zhukov1, zhukov2}. In order to formulate a precise question,
it may be easiest to pass to the context of 
considering a representation of the \'etale
fundamental group of a surface and computing its conductor along different
boundary divisors. Indeed, this will be the point of view of the sequel to
this paper.

There is also a construction of Artin conductors in the imperfect
residue field case due to Borger \cite{borger}, by passing from
$E$ to a certain extension which is universal for the property of
having perfect residue field. Borger's construction does not behave well
with respect to tame base extension, but one should get a better invariant
by forcing such good behavior (i.e., constructing a logarithmic
analogue of Borger's conductor). Indeed, we expect the
following.

\begin{conj} \label{conj:borger}
For $\rho$ a representation with finite local monodromy, for $m$ a
positive integer, let $E_m$ be an extension of $E$ 
which is tamely ramified of tame degree $m$.
Let $b'(E_m)$ be Borger's Artin
conductor of the restriction of $\rho$ to $G_{E_m}$.
Then the limsup of $m^{-1} b'(E_m)$ as $m \to \infty$ equals
the differential Swan conductor of $\rho$.
\end{conj}

Since the Abbes-Saito construction also works in mixed characteristic,
one may also be interested in reconciling it there with a differential
construction. For more on this possibility, see the next subsection.

\subsection{Comments on mixed characteristic}
\label{subsec:mixed}

It would be interesting to extend the constructions in this paper to
the case where $R$ has mixed characteristics. The analogue of the passage from
Galois representations to $\nabla$-modules is given by $p$-adic Hodge theory,
specifically via the theory of $(\phi, \Gamma)$-modules over the Robba ring, 
as in the work of Fontaine, Cherbonnier-Colmez, Berger, et al.

In that context, when $R$ has perfect residue field,
Colmez \cite{colmez} has given a recipe for reading
off the Swan conductor of de Rham representations from the
associated $(\phi, \Gamma)$-module.
One would like to reformulate this recipe via Berger's
construction of the Weil-Deligne representation, which converts
the $(\phi, \Gamma)$-module into a $(\phi, \nabla)$-module over $\calR$
\cite{berger}; however, it is not immediately clear how to do this.
The fact that this might even be possible is suggested by work of
Marmora \cite{marmora}, who gives a direct comparison
with differential Swan conductors, but only for the Swan conductor of
a representation over the maximal $p$-cyclotomic extension of a given
$p$-adic field.

If one can indeed give a differential definition of the usual Swan 
conductor of a de Rham representation in the perfect residue field,
then it seems likely one can make a differential definition in the
imperfect residue field case. Indeed, the construction of 
$(\Phi, \Gamma)$-modules has already been generalized to this setting
by Morita \cite{morita}.
If one can do all this, 
then one will again encounter the question of reconciliation with the
Abbes-Saito constructions; however, it is not clear whether in this case
the Hasse-Arf theorem would be any easier on the differential side than
on the Abbes-Saito side.

\end{document}